\documentclass[11pt,twoside]{article}
\usepackage{latexsym, amssymb, amsmath, theorem}
\usepackage{epsfig}
\numberwithin{equation}{section}
                                                                           
\theoremstyle{plain}
\theorembodyfont{\itshape}
\newtheorem{theorem}{Theorem}[section]
\newtheorem{proposition}[theorem]{Proposition}
\newtheorem{lemma}[theorem]{Lemma}
\newtheorem{corollary}[theorem]{Corollary}
                                                                                
\theorembodyfont{\rmfamily}
\newtheorem{definition}[theorem]{Definition}

\newtheorem{example}[theorem]{Example}
\newtheorem{remark}[theorem]{Remark}
\newenvironment{proof}{{\noindent \textbf{Proof}\,\,}}
{\hspace*{\fill}$\Box$\medskip}

\title{On density of horospheres in dynamical laminations
\protect\footnote{AMS Classification 58F23, 57M50}}
                                                                                
\author{A. Glutsyuk
\thanks{Poncelet Laboratory (UMI of CNRS and Independent University of Moscow), 11, 
Bolshoi Vlasievskii Pereulok, 119002, Moscow, Russia} \thanks{Permanent address: 
CNRS, Unit\'e de Math\'ematiques Pures et Appliqu\'ees,
M.R., \'Ecole Normale Sup\'erieure de Lyon, 46 all\'ee d'Italie,
69364 Lyon 07, France. \newline Email:
aglutsyu@umpa.ens-lyon.fr}
\thanks{Supported by part by RFBR grants 02-02-00482, 
07-01-00017-a 
and NTsNIL\_a (RFBR-CNRS)  05-01-02801 and by ANR grant ANR-08-JCJC-0130-01.}
}

\begin{document}
\maketitle
\def\hw{\hat w}                                                                             
\def\Cal{\mathcal}
\def\Bbb{\mathbb}
\def\e {\varepsilon }
\def\nbd {neighborhood }
\def\nbds {neighborhoods }
\def\ph {\varphi }
\def\a {\alpha }
\def\de {\delta }
\def\g {\gamma }
\def\rr{\mathbb R}
\def\ca{\Cal A}
\def\caf{\Cal A_f}
\def\can{\caf^n}
\def\ch{\Cal H}
\def\chf{\ch_f}
\def\cc{\mathbb C}
\def\oc{\overline\cc}
\def\rr{\mathbb R}
\def\hf{\hat f}
\def\chff{\chf\slash\hf}
\def\ha{\hat a}
\def\hz{\hat z}
\def\hb{\hat b}
\def\hx{\hat x}
\def\hh{\mathbb H^3}
\def\zz{\mathbb Z}
\def\var{\varepsilon}
\def\b{\beta}
\def\th{\tilde h}
\def\nn{\mathbb N}
\def\htr{\mathbb H^3}
\def\hc{\hat c}
\def\hy{\hat y}
\def\hw{\hat w}
\def\chffp{\chf'\slash\hf}

\tableofcontents
                                                                          
\section{Introduction and main results}

\subsection{Introduction, brief description of main results and the plan of the paper}

In 1985 D.Sullivan \cite{su} had introduced a dictionary between two domains of 
complex dynamics: iterations of rational functions 
$f(z)=\frac{P(z)}{Q(z)}:\oc\to\oc$ on the Riemann sphere and Kleinian groups. 
The latter are discrete subgroups of the group of conformal automorphisms 
of the Riemann sphere.  This dictionary motivated many remarkable results in both domains, 
starting from the famous Sullivan's no wandering domain theorem \cite{su} in the theory of iterations  
of rational functions. 

One of the principal objects used in the study of Kleinian groups is the  hyperbolic 
3- manifold associated to a Kleinian group, which is the quotient of its lifted action to the 
hyperbolic 3- space $\htr$. 
M.Lyubich and Y.Minsky have suggested to extend Sullivan's dictionary by 
providing an analogous construction for iterations of 
rational functions. For each rational function $f$  they have 
constructed a {\it hyperbolic lamination} \ $\chf$ 
(see \cite{lm} and Subsection 1.3 below). This is a topological space 
foliated by hyperbolic 3- manifolds (some of them may have singularities) so that 

- a neighborhood of a  nonsingular point is fiberwise homeomorphic to the product of 
a part of the Cantor set and 3- ball;

- the  hyperbolic metric of leaves depends continuously on the transverse parameter;

- there exists a natural projection $\chf\to\oc$ under which the (non-bijective) 
action $f:\oc\to\oc$ lifts up to a homeomorphic action 
$\hf:\chf\to\chf$ that maps leaves to leaves isometrically;

- the lifted action $\hf$ is proper discontinuous, and hence, the quotient 
$\chff$ is a Hausdorff topological space laminated by hyperbolic 3- manifolds 
(called {\it the quotient hyperbolic lamination});

- the lamination $\chf$ either is minimal itself (i.e., each leaf is dense), or 
becomes minimal after removing a finite number of isolated leaves 
(which may exist only in very exceptional cases); the complement to the 
isolated leaves will be denoted $\chf'$. 

The hyperbolic lamination $\chf$ is constructed as follows. Take the natural extension 
of the dynamics of $f$ to the space of all its backward orbits. The latter space always 
contains some Riemann surfaces conformally-equivalent to $\cc$. The union of all these 
surfaces (denoted by $\caf^n$) is invariant under the lifted dynamics. Pasting a copy 
of the hyperbolic 3- space to each surface, appropriate strengthening of the topology 
and completion of the new space thus obtained yields the hyperbolic lamination $\chf$. 
 
 Recent studies of the hyperbolic 3-manifolds associated to Kleinian groups 
 resulted in resolving of all big problems of the theory, including a 
 positive solution of the famous Ahlfors measure conjecture  
 (with contributions of many mathematicians, see the papers 
\cite{ah1}, \cite{ah2} and references therein).  
On the other hand, very recently the analogous conjecture in the 
theory of rational iterations was proved to be wrong: 
by using a completely different idea proposed by  A.Douady,  
X.Buff and A.Cheritat have constructed 
examples of quadratic polynomials with Julia sets of positive measure \cite{bc}.

There is a hope that similarly,   
hyperbolic laminations associated to rational functions would shed new light 
on the underlying dynamics. 

One of the main problems in the theory of rational iterations is the Fatou conjecture: 
{\it is it true that each critical orbit of a generic rational function of a given degree 
either is periodic itself, or 
converges to an attracting or super-attracting periodic orbit?} 
It is open already for the quadratic polynomials, for which it  
is equivalent to the No Invariant Line Field Conjecture: {\it is it true that there are no measurable invariant line fields supported on the Julia set?}

The answer to the latter conjecture is known to be positive for the 
critically-nonrecurrent rational functions  
without parabolics that are different from the Latt\`es examples.    
Elementary geometric proofs using laminations may be found in   
\cite{lm} (proposition 8.9) and \cite{g2} (subsection 4.3). 

There is a hope that the laminations could be helpful in studying the No 
Invariant Line Field Conjecture. 

The present paper studies the arrangement of the horospheres in the quotient hyperbolic 
lamination $\chff$. Let us recall their definition. 
The hyperbolic space with a marked point ``infinity'' on its boundary Riemann sphere 
 admits a standard model of  
half-space  in the Euclidean 3- space. Its isometries that fix the infinity are exactly the extensions 
of the complex affine transformations of the boundary (we call these extensions 
``affine isometries'').  
A horizontal plane (i.e., a plane parallel to the boundary) in the  half-space is called a {\it horosphere}.  The affine isometries transform the 
horospheres to the horospheres. The {\it horospheres of a quotient of} $\htr$ by a discrete 
group of affine isometries are the quotient projection images 
 of the horospheres in $\htr$. All the  horospheres mentioned above carry natural complex 
 affine structures (they may have conical singularities) and foliate the 
 ambient hyperbolic manifold (orbifold). 

The leaves of the laminations $\chf$ and $\chff$ are also quotients of 
$\htr$ by discrete groups of affine isometries. Thus, all their  leaves are 
foliated by well-defined horospheres. 

The arrangement of the horospheres in the above hyperbolic laminations is related 
to the behavior of the modules of the derivatives of the iterations of the rational function. 

The vertical geodesic flow acts on $\htr$ by translating points along the 
geodesics issued from the infinity. This yields the leafwise vertical geodesic 
flows acting on $\chf$ and $\chff$, for which the horospheric laminations are 
the unstable laminations. 

The classical results concerning the geodesic flows on compact hyperbolic 
surfaces say that the horocyclic lamination is minimal \cite{he} and 
uniquely ergodic \cite{furst, m}. Their generalizations have found important applications 
in different domains of mathematics, including number theory. We hope that 
studying the vertical geodesic flow and the horospheric lamination would 
have applications in understanding the underlying dynamics. 

The main results are stated in 1.5 and proved in Sections 2-5.  
The principal one (Theorem \ref{trep} proved in Section 2) 
says that at least some horospheres are always 
dense in $\chffp$, provided that the mapping $f$ does not belong to the 
following list of exceptions: 
\begin{equation}z^{\pm d}, \ \text{Chebyshev polynomials, \ Latt\`es examples.}\label{list}\end{equation}
Namely, all the horospheres ``over'' the so-called 
branch-nonexceptional repelling periodic orbits are always dense. 

\begin{remark} For any $f$ from the list (\ref{list}) each horosphere in $\chff$ is nowhere 
dense. This follows from a result of \cite{ly} (see Corollary \ref{nodense} in 1.5). 
\end{remark}

Theorem \ref{thyp} (proved in 1.5)  says that all the horospheres 
are dense in $\chffp$, if $f$ does not belong to the list (\ref{list}) and is 
critically-nonrecurrent without parabolic periodic points (e.g., 
hyperbolic). In the case, when $f$ does not belong 
to (\ref{list}),  is critically-nonrecurrent 
and has parabolic periodic points, Theorem \ref{tnrec} (also proved in 1.5) 
says that all the horospheres are dense in $\chffp$, 
except for the horospheres ``related'' to the parabolic points. To prove density of  a horosphere in Theorem \ref{tnrec}, 
we show (Theorem \ref{taccum} stated in 1.5 and proved in Section 3) that it  
accumulates to some horosphere over appropriate branch-nonexceptional 
repelling periodic orbit.  
(The limit horosphere is dense by Theorem \ref{trep}.) Theorem \ref{taccum2} (stated in 1.5 and 
proved in Section 4) deals with an arbitrary rational function having 
a parabolic periodic 
point. It says that each horosphere in a leaf associated to this point is closed in $\chff$ 
and does not accumulate to itself. 

\begin{remark} While the present paper was in preparation, Theorem \ref{thyp} has 
already been applied to prove the unique ergodicity of the  quotient 
horospheric laminations associated to appropriate rational functions. The latter 
include  all the hyperbolic and critically-finite functions \cite{gll}. 
The results of the paper \cite{gll} with brief proofs are announced in \cite{g2}.  
\end{remark}

\begin{remark} There exist (even hyperbolic) rational functions that do not belong 
to the list (\ref{list}) such that the corresponding 
hyperbolic lamination $\chf$ has a leaf whose  horospheres are nowhere dense in 
$\chf$. This is true, e.g., for the quadratic polynomials 
$f_{\var}(z)=z^2+\var$ with $\var\in E=(-\infty,\frac14)\setminus\{0,\ -2\}$, 
and also with $\var$ belonging to a complex neighborhood of the set $E$ 
(Theorem \ref{tnodense} and its Addendum, both stated in 1.5 and proved in \
Section 5). The above leaf 
with nondense horospheres is associated to a repelling fixed point 
(that is real, if so is $\var$). 
\end{remark} 

\begin{example} Consider the above quadratic polynomial family 
$f_{\var}(z)=z^2+\var$. It is well-known that the quotient laminations 
$\Cal H_{f_0}\slash\hf_0$ and $\Cal H_{f_{\var}}\slash\hf_{\var}$ are homeomorphic 
for all $\var\neq0$ small enough. 
The homeomorphism sends leaves to leaves but not isometrically. 
On the other hand, Theorem \ref{thyp} implies that if $\var\neq0$ is small enough, 
then each horosphere in the latter lamination is dense, while no horosphere 
in the former one is dense (Corollary \ref{nodense}). 
\end{example}

The necessary background material is recalled in Subsections 1.2 
(rational iterations), 
1.3 (affine and hyperbolic laminations) and 1.4 (horospheres and their metric properties). 

For the proof of Theorem \ref{trep} we fix a horosphere in $\chf$ ``over'' a 
branch-nonexceptional repelling periodic orbit  and show that the union of the images of the horosphere 
under the forward and the backward iterations of $\hf$ 
 is dense in $\chffp$. To do this, we study the holonomies of the horosphere along loops 
based at a repelling periodic point. We show that  the images of a point of the horosphere under subsequently applied dynamics 
and holonomies are dense in the projection preimage of the base point. To this 
end, we use the description of  the holonomy in terms of  
basic cocycle (its definition and some basic properties are recalled in Subsection 1.4). 

Some results of the paper with brief proofs were announced in \cite{gl}. 

Earlier some partial result on density of horospheres was obtained in a joint work by 
M.Yu.Lyubich and D.Saric \cite{lypr} (under additional assumptions on 
the arithmetic nature of the multipliers at the repelling periodic points). 

Everywhere below we assume that the rational function $f=\frac{P(z)}{Q(z)}:\oc\to\oc$ 
under consideration  has degree at least 2.

\def\nn{\mathbb N}
\def\crf{\Cal R_f}

\subsection{Background material 1: rational iterations} 
The basic notions and facts of holomorphic dynamics 
recalled here are contained, e.g., in \cite{lyu} and \cite{lm}. Let 
$$f=\frac{P(z)}{Q(z)}:\oc\to\oc \ \text{be a rational function. Recall that}$$ 

- its {\it Julia set} $J=J(f)$ is the closure of the union of  
the repelling periodic points, see  the next definition. An equivalent definition of 
the Julia set says that its complement $\oc\setminus J$ (called the {\it Fatou set}) 
is the maximal open subset where the iterations $f^n$ form a normal family 
(i.e., are equicontinuous on compact subsets). One has 
\begin{equation}f^{-1}(J)=J=f(J).\label{julin}\end{equation}

\begin{definition} A germ of nonconstant 
holomorphic mapping $f:(\cc,0)\to(\cc,0)$ at a fixed point 0 is called 
{\it attracting (repelling / parabolic, superattracting)}, if its derivative at the fixed point 
respectively has nonzero module less than 1 
(has module greater than 1 / is equal to a root of unity and no iteration of the 
mapping $f$ is identity / is equal to zero). An attracting (repelling, parabolic or 
superattracting) periodic point of a rational mapping is 
a fixed point (of the corresponding type) of its iteration. 
\end{definition}

\begin{definition} A rational function is said to be {\it hyperbolic}, if the forward orbit of each its critical point either is periodic itself (and hence, superattracting), or tends to an attracting (or a superattracting) periodic orbit.\end{definition}

\begin{definition} Given  a rational function. A point of the Riemann sphere is called {\it postcritical}, if  
it belongs to the forward orbit of  a critical point. A rational function is called 
{\it critically-finite}, if the number of its postcritical points is finite. 
\end{definition}

\def\o{\omega}

\begin{definition} The $\o$- limit set $\o(c)$ of a point $c\in\oc$ is the set of limits of 
converging subsequences of its forward 
orbit $\{ f^n(c) | n\geq0\}$ (the $\omega$- limit set of a periodic orbit is the orbit itself). 
A point $c$ is called {\it recurrent}, if $c\in\o(c)$. 
\end{definition}

\begin{definition} A rational mapping is called {\it critically-nonrecurrent}, if each its critical point is 
either nonrecurrent, or periodic (or equivalently, each critical point in the Julia set is nonrecurrent). 
\end{definition}

\begin{example} The following mappings are critically-nonrecurrent: any hyperbolic mapping; 
any critically-finite mapping; any quadratic polynomial with a parabolic periodic orbit. 
A hyperbolic mapping has no parabolic periodic points. 
\end{example}

\begin{theorem} (\cite{lyu}, p.60) 
A germ of conformal mapping at an attracting (repelling) fixed point is always 
conformally linearizable: there exists a local conformal coordinate in which the germ is equal 
to its  linear part (the multiplication by its derivative at the fixed point). 
\end{theorem}

\begin{remark} \label{parab} 
Let $f(z)=z+z^{k+1}+\dots$ be a parabolic germ tangent to the identity. The set $\{ z^k\in\rr_+\}$ 
consists of $k$ rays going out of 0 (called {\it repelling rays}) such that 

- each repelling ray is contained in appropriate sector $S$ (called {\it repelling sector}) for which 
there exists an arbitrarily small  neighborhood $U=U(0)\subset\cc$ where $f$ is univalent and 
such that $f(S\cap U)\supset S\cap U$  
and each backward orbit of the restriction $f|_{S\cap U}$  enters the fixed point 0 asymptotically 
along the corresponding repelling ray;

- there is a canonical 1-to-1 conformal coordinate $t$ on $S\cap U$ in which 
$f$ acts by translation: $t\mapsto t+1$; if the above sector $S$ is chosen large enough, 
then this coordinate parametrizes $S\cap U$ 
by a domain in $\cc$ containing a left half-plane; the above coordinate is  well-defined 
up to translation and is called {\it Fatou coordinate} (see \cite{ec}, \cite{vor}). 

For any parabolic germ (not necessarily tangent to the identity) its appropriate iteration is tangent 
to the identity. By definition, the repelling rays and sectors of the former are those (defined above) 
of the latter. 
\end{remark}

Let us recall what are Chebyshev polynomials and Latt\`es examples. 

{\bf Chebyshev polynomials.} For any $n\in\nn$ there exists a unique (real) 
polynomial $p_n$ of degree $n$ 
that satisfies the trigonometric identity $\cos n\theta=p_n(\cos\theta)$. It is called {\it Chebyshev 
polynomial}.

{\bf Latt\`es examples.} Consider a one-dimensional 
complex torus, which is the quotient of $\cc$ by a lattice. Consider arbitrary 
multiplication by a constant $\lambda\in\cc$, $|\lambda|>1$, 
that maps the lattice to itself. It induces an endomorphism of the torus of degree greater 
than 1. The quotient of the torus by the central symmetry $z\mapsto -z$ is 
the Riemann 
sphere.  The above endomorphism together with the quotient projection induce a 
rational transformation of the Riemann sphere called {\it Latt\`es example}. 

\begin{remark} Let $f$ be either Chebyshev, or Latt\`es. Then it is critically finite. More precisely, 
the forward critical orbits eventually finish at repelling fixed points. The Julia set of a Chebyshev polynomial is the segment $[-1,1]$ of the real line, while that of a Latt\`es example is the whole Riemann sphere.
Chebyshev and Latt\`es functions have branch-exceptional repelling fixed points, 
see the following definition. 
\end{remark}

\begin{definition} \label{crcom} \cite{klr} A periodic point of a 
rational function 
is called {\it branch-exceptional}, if any its nonperiodic backward orbit contains a 
critical point. In this case its periodic orbit is also called 
branch-exceptional.  
\end{definition}

\begin{remark} (Lasse Rempe \cite{klr}). There exist quadratic rational 
functions with a branch-exceptional repelling fixed point that are neither 
Chebyshev, nor Latt\`es.
\end{remark}

\subsection{Background material 2: affine and hyperbolic dynamical laminations}
The constructions presented here were introduced in \cite{lm}. We recall them briefly and 
send the reader to \cite{lm} for more details. 

Recall that a {\it lamination} (by manifolds) is a ``topological'' foliation by manifolds, i.e., a 
topological space that is split as a disjoint union of manifolds (called {\it leaves}) 
of one and the same dimension so that each point of the ambient space admits 
a neighborhood (called ``flow-box'') such that each connected component 
(local leaf) of its intersection 
with each leaf is  homeomorphic to a ball; the neighborhood itself is 
homeomorphic to the product of the ball  and some (transverse) topological space under a 
 homeomorphism transforming the local leaves to the fibers of the product. 

Let $f:\oc\to\oc$ be a rational function. Set  
$$\Cal N_f=\{\hz=(z_0,z_{-1},\dots) \ | \ z_{-j}\in\oc, \ f(z_{-j-1})=z_{-j}\}. $$
This is a topological space equipped with the natural product topology and the projections 
$$\pi_{-j}:\Cal N_f\to\oc, \ \hz\mapsto z_{-j}.$$
The action of $f$ on the Riemann sphere lifts naturally up to a homeomorphism 
$$\hf:\Cal N_f\to \Cal N_f, \ (z_0,z_{-1},\dots)\mapsto(f(z_0),z_0,z_{-1},\dots), \ 
f\circ\pi_{-j}=\pi_{-j}\circ\hf.$$
First of all we recall the construction of  the ``regular leaf subspace'' 
$\Cal R_f\subset\Cal N_f$, which is 
a union of Riemann surfaces that foliate $\Cal R_f$ in a very turbulent way. Afterwards we take 
the subset $\caf^n\subset \Cal R_f$ of the leaves conformally-equivalent to $\cc$. Then we refine 
the induced topology on $\caf^n$ to make it a lamination (denoted $\caf^l$) by complex lines with 
a continuous family of affine structures on them. Afterwards we take a completion  
$\caf=\overline{\caf^l}$ in the new topology. The space $\caf$ is a lamination by affine 
Riemann surfaces (the new leaves added by the completion may have conical singularities). 
Then we discuss 
the three-dimensional extension of $\caf$ up to a lamination $\chf$ by hyperbolic manifolds 
(with singularities). 

Let $\hz\in\Cal N_f$, $V=V(z_0)\subset\oc$ be a neighborhood of $z_0$. For any 
$j\geq0$ set 
$$V_{-j}=\text{the connected component of the preimage} \ f^{-j}(V) \ 
\text{that contains} \ z_{-j}.$$
$$ \text{Then} \ V_0=V, \ \text{and} \ f^j:V_{-j}\to V \ \text{are ramified coverings}.$$
 
\begin{definition} We say that a point $\hz\in\Cal N_f$ is {\it regular}, if 
there exists a disk $V$ containing 
$z_0$ such that the above coverings $f^j:V_{-j}\to V$ have uniformly bounded 
degrees. Set
$$\Cal R_f\subset\Cal N_f \  \ \text{the set of the regular points in} \ \Cal N_f.$$
\end{definition}

\begin{example} Let $\hz\in\Cal N_f$ be a backward orbit such that there exists a 
$j\in\nn\cup0$ for which  the point $z_{-j}$ is disjoint from the $\o$- limit sets of the critical points. 
Then  $\hz\in\crf$. If the mapping $f$ is hyperbolic, then this is the case, if and only if $\hz$ 
is not a (super) attracting periodic orbit. A mapping $f$ is critically-nonrecurrent, if and only if 
$$\crf=\Cal N_f\setminus\{\text{attracting and parabolic periodic orbits}\}, \ \text{see \cite{lm}.}$$
\end{example}

\begin{definition} Let $\hz\in\crf$, $V$, $V_{-j}$ be as in the above definition. 
The {\it local leaf} $L(\hz,V)\subset\crf$  is the set of the 
points $\hz'\in\crf$ such that $z'_{-j}\in V_{-j}$ for all $j$ (the local leaf is path-connected by definition). 
We say that the above local leaf is {\it univalent over} \ $V$, if the projection 
$\pi_0$ maps it bijectively onto $V$. 
The {\it global leaf}  containing $\hz$  
(denoted $L(\hz)$)  is the maximal path-connected subset in $\crf$ containing $\hz$. 
\end{definition}

\begin{remark} Each leaf $L(\hz)\subset\crf$ carries a natural structure 
of  Riemann surface so that  the  restrictions to the leaves of the above projections $\pi_{-j}$ 
are meromorphic functions. A local leaf $L(\hz,V)\subset L(\hz)$ 
(when well-defined) is the connected  component containing $\hz$ of the preimage 
$(\pi_0|_{L(\hz)})^{-1} (V)\subset L(\hz)$. 
\end{remark}

\begin{remark} \label{itersame} 
The above-defined objects $\crf$, $\Cal R_{f^n}$ corresponding to both $f$ and any its 
forward iteration $f^n$, are naturally homeomorphic under the mapping that sends a 
backward orbit $\hz\in\Cal N_f$ to the backward orbit $(z_0,z_{-n},z_{-2n},\dots)
\in\mathcal N_{f^n}$. 
The latter homeomorphism maps the leaves conformally 
onto the leaves. 
\end{remark}

Recall that given a rational function $f$ and a $n\in\mathbb N$, 
a $n$- periodic connected component $U$ of the Fatou set is called 
a {\it rotation domain: Siegel disk (Herman ring)}, if the restriction 
$f^n|_U$ is conformally conjugated to a rotation of disk (annulus). 

We use the following 
\begin{lemma} \label{shr} {\bf (Shrinking Lemma, \cite{lm}, p.86)} Let $f$ be a rational mapping, $V\subset\oc$ 
be a domain, $V'\Subset V$ be a compact subset. Then for any sequence of 
single-valued inverse branches 
$f^{-n}:V\to\oc$  the diameters of the images $f^{-n}(V')$ tend to 0, as 
$n\to+\infty$, if ond only if $f$ has no rotation domain 
that contains an infinite number of the above images.
\end{lemma}

\begin{definition} A noncompact Riemann surface is said to be 
{\it hyperbolic (parabolic),} 
if its universal covering is conformally equivalent to the unit disk (respectively, $\cc$).
\end{definition}
\begin{remark} Parabolic leaves in $\crf$ always exist (see the next two 
examples) and are simply connected; 
hence they  are conformally equivalent to $\cc$ (\cite{lm}, 
lemma 3.3, p.27). If $f$ is critically-nonrecurrent, then 
each leaf is parabolic (\cite{lm}, proposition 4.5, p.36). 
On the other hand, there are rational 
mappings such that some leaves of $\crf$ are 
hyperbolic (e.g., the mappings with rotation domains, see \cite{lm}, p.27). 
Nontrivial 
examples of rational functions without rotation domains and with infinitely 
many hyperbolic leaves in $\crf$, whose projections intersect the Julia set,  
were constructed in \cite{klr}. 
\end{remark}

\begin{example} \label{exrep} 
 Let $a\in\oc$ be a repelling fixed point of $f$, $\ha=(a,a,\dots)\in\Cal N_f$ 
be its fixed orbit. Then $\ha\in\crf$ and the leaf $L(\ha)$ is parabolic (it is $\hf$- invariant and the 
quotient of $L(\ha)\setminus\ha$ by $\hf$ is a torus).
The linearizing coordinate 
$w$ of $f$ in a neighborhood of $a$ lifts up to a conformal isomorphism 
$w\circ\pi_0:L(\ha)\to\cc.$ 
Analogously, the periodic orbit of a repelling periodic 
point is contained in a parabolic leaf (see  Remark \ref{itersame}).  
\end{example}

\begin{example} Let $f$ have a parabolic fixed point $a\in\oc$, $f'(a)=1$, 
$\ha=(a,a,\dots)\in\Cal N_f$ be its 
fixed orbit. Then $\ha\notin\crf$. On the other hand, for each repelling ray 
(see Remark \ref{parab}) there is a unique leaf in $\crf$ (denoted $L_a$) 
consisting of the backward orbits that converge to $a$ asymptotically along the 
chosen ray. The leaf $L_a$ is parabolic: the Fatou coordinate $w$ on the 
corresponding repelling sector lifts up to a conformal isomorphism 
$w\circ\pi_0:L_a\to\cc$. 
An analogous statement holds true in the case, when $a$ is  a 
parabolic periodic point and $f$ is not necessarily tangent to the identity 
there. 
\end{example}
\begin{definition} \label{defass1} The leaves from the two above examples are  called respectively a 
{\it leaf associated to a repelling (parabolic) periodic point.}
\end{definition}
\begin{proposition} \label{repell} 
 A point $\hz\in\Cal N_f$ belongs to a leaf associated to a repelling 
 (or parabolic) fixed point $a$,  if and only if it is represented by a 
 backward orbit converging to $a$ 
 (and distinct from its fixed orbit, if the latter is parabolic). 
\end{proposition}

The proposition follows from  the Shrinking Lemma. 

\begin{corollary} \label{coruniq} A leaf of $\Cal N_f$  
 can be associated to at most one repelling or parabolic periodic point.
 \end{corollary}

Set
$$\can=\ \text{the union of the parabolic leaves in} \ \crf.$$
If $f$ is hyperbolic, then $\caf^n$ is a  lamination with a global Cantor 
transverse section.  In general, $\caf^n$ is not a lamination in a good sense, 
since some ramified local leaves can 
accumulate to a univalent one in the product topology. The refined topology (defined in \cite{lm}) 
that makes it an orbifold lamination is recalled below. To do this, we use the following 

\begin{remark} Let $\hz\in\caf^n$. Fix a conformal isomorphism $\cc\to L(\hz)$ that 
sends 0 to $\hz$ (it is unique up to multiplication by nonzero complex constant in 
the source). For any $j\geq0$ set 
\begin{equation}\phi_{-j,\hz}=\pi_{-j}|_{L(\hz)}:\cc=L(\hz)\to\oc; \ 
\phi_{-j,\hz}=f\circ\phi_{-j-1} \ \text{for any} \ j, \ 
\phi_{-j,\hz}(0)=z_{-j}.
\label{allfunc}\end{equation}
This is a  meromorphic function sequence, uniquely defined up to the 
$\cc^*$- action on the source space $\cc$ by multiplication by complex 
constants.  Two points of $\caf^n$ lie in one and the same 
leaf, if and only if the corresponding function sequences are obtained from each other 
by affine transformation of the source variable. 
\end{remark}

\def\chk{\hat{\Cal K}_f}
Let $\chk$ denote the space  of the nonconstant 
meromorphic function sequences 
\begin{equation}\{\phi_{-j}(t)\}_{j\in\nn\cup0}, \ \phi_{-j}:\cc\to\oc,\ \phi_{-j}=f\circ\phi_{-j-1} \ \text{for all}\ j.
\label{funcsec}\end{equation}
This is a subset of the infinite product of copies of the meromorphic function space; the latter space 
is equipped with the topology of uniform convergence on compact sets. The product topology induces a topology on the space $\chk$. The groups $Aff(\cc)$ 
(complex affine transformations of $\cc$), 
$\cc^*\subset Aff(\cc)$ and $S^1=\{|z|=1\}\subset\cc^*$ act on the space 
$\chk$  by variable changes in the source. The stabilizer in $Aff(\cc)$ 
of a meromorphic function sequence from $\chk$ 
is a discrete group of Euclidean isometries. Set 
\begin{equation}\chk^a=\chk\slash\cc^*, \ \chk^h=\chk\slash S^1.
\label{graction}\end{equation}
The spaces $\chk^a$, $\chk^h$ are equipped with the quotient topologies 
induced from $\chk$. A {\it leaf} in $\chk^a$ ($\chk^h$) is the quotient 
projection of an $Aff(\cc)$- orbit in $\chk$. 

\begin{remark} \label{afcont} A leaf in $\chk^a$ ($\chk^h)$ is 
 naturally identified  
with the quotient of $Aff(\cc)\slash\cc^*=\cc$ (respectively, 
$Aff(\cc)\slash S^1=\htr$) 
by the left action of a discrete subgroup of Euclidean isometries in $Aff(\cc)$. 
This equips the leaves 
with natural {\it affine (hyperbolic) structures} that {\it vary continuously} 
on the ambient space.  Each leaf is thus an affine (hyperbolic) orbifold. 
There is  a natural inclusion 
\begin{equation}\caf^n\to\chk^a: \ \hz\mapsto\{\phi_{-j,\hz}\}_{j=0}^{+\infty}\slash\cc^*, 
\text{ see (\ref{allfunc}).}\label{canincl}\end{equation}
\end{remark}
The action of $\cc^*=S^1\times \rr_+^*$ on $\chk$ induces an action of $\rr_+^*$  on 
$\chk^h$. The quotient of the latter action is 
$\chk^a=\chk\slash\cc^*=\chk^h\slash\rr_+^*$. The corresponding quotient projection 
will be denoted
 \begin{equation}\pi_{hyp}:\chk^h\to\chk^a. \label{pih}\end{equation}
 The projection $\pi_{hyp}$ maps each hyperbolic leaf of $\chk^h$ onto an 
 affine leaf of $\chk^a$ that is canonically identified with the boundary of the hyperbolic 
 leaf. Conversely, the preimage of an affine leaf is a hyperbolic leaf.
 
\begin{definition} \label{deflam}
The topological subspace $\caf^l\subset\chk^a$ is the image of the space $\caf^n$ under the inclusion (\ref{canincl}), or equivalently, the space $\caf^n$ equipped with 
the topology induced from $\chk^a$. The space $\caf$ is the closure of 
$\caf^l$ in the space $\chk^a$. We set  
\begin{equation}\chf^l=\pi_{hyp}^{-1}(\caf^l)\subset\chk^h, \ \chf=\pi_{hyp}^{-1}(\caf)=
\overline{\chf^l}\subset\chk^h.\label{defch}\end{equation}
The space $\caf$ ($\chf$) is called the {\it affine (hyperbolic) orbifold lamination 
associated to} $f$. For any point $\hz\in\caf$ the affine leaf through $\hz$ in $\caf$ 
will be denoted by $L(\hz)$, and set  
$$H(\hz)=\pi_{hyp}^{-1}(L(\hz))\subset\chf: \text{ the corresponding hyperbolic leaf.}$$
\end{definition}

\begin{remark} \label{remquot} 
In general, the topology of the space $\caf^l$ is stronger than that of $\caf^n$.  
The spaces $\caf^l$, $\caf$, $\chf^l$, $\chf$ consist of entire leaves. Each leaf 
in $\caf^l$ ($\chf^l$) is affine-isomorphic (isometric) to $\cc$ (respectively, $\htr$). 
Other leaves may contain conical singularities. 
\end{remark}

There is  a natural  projection 
\begin{equation}p:\chk^a\to\caf^n\label{prp}\end{equation} 
induced by the mapping $\chk\to\caf^n$ that sends each 
sequence (\ref{funcsec}) of functions to the sequence of their values at 0. 
The latter 
value sequence is always a regular backward orbit of $f$, and it lies  in a parabolic leaf of 
$\crf$. The regularity follows from definition. The parabolicity  
follows from Picard's theorem. The composition of $p$ with the 
natural inclusion $\caf^n\to\chk^a$ is the identical mapping $\caf^n\to\caf^n$. 
The projection 
\begin{equation}\caf\to\oc \ \text{induced 
by}  \ \pi_0, \ (\phi_{-j})_{j\in\nn\cup0}\mapsto\phi_0(0), \ \text{will be also denoted by}  \ \pi_0.\label{pi0ext}
\end{equation}
  The rational mapping  $f:\oc\to\oc$ lifts up to the leafwise homeomorphism 
 $$\hf:\chk\to\chk,  \ \hf:(\phi_0,\phi_{-1},\dots)\mapsto(f\circ\phi_0,\phi_0,\phi_{-1},\dots), \ \text{which induces homeomorphisms} $$
 $$\hf:\caf\to\caf \ \text{leafwise affine \ and} \ \hf:\chf\to\chf \ 
 \text{leafwise isometric; } \hf(\caf^l)=\caf^l.$$
  These homeomorphisms form a commutative diagram with the projection $\pi_{hyp}$. 
  The latter homeomorphism acts  properly discontinuously on $\chf$, and its quotient 
 $$\chff \ \text{is Hausdorff and 
 called the {\it quotient hyperbolic lamination associated to}}\ f.$$
 
  \begin{proposition}  \label{conva} 
(\cite{lm}, proposition 7.5, p.62) 
A sequence of points $\hz^m\in\caf^l$ converges to a point $\hz\in\caf^l$, as $m\to\infty$, 
if and only if 

 -  \ $\pi_{-j}(\hz^m)\to\pi_{-j}(\hz)$ for any $j$, 

 -  \  for every $N\in\mathbb N$, each connected domain $V\subset\oc$ and  
 each its subdomain $U$ such that  $\overline U\subset V$, $\pi_{-N}(\hz)\in U$, 
   if the local leaf $L(\hf^{-N}(\hz),V)$ is univalent over $V$, then the local leaf 
 $L(\hf^{-N}(\hz^m),U)$ is univalent over $U$, whenever $m$ is large enough. 
 \end{proposition}
 
  \begin{remark} \label{remcaf} 
  The analogous criterion holds true for convergence of a sequence of points 
 in $\caf$ to a point in $\caf^l$ with  the following definition of local leaf in $\caf$. 
 \end{remark}
 
 \begin{definition} \label{univdefn} 
 Let $f$ be a rational mapping, $\caf$ be the corresponding affine lamination, 
 $L\subset\caf$ be a leaf, $\hz\in L$, $V\subset\oc$ be a domain containing its projection $\pi_0(\hz)$. 
 The {\it local leaf} $L(\hz,V)$ is the connected component containing $\hz$ of the projection preimage 
 $\pi_0^{-1}(V)\cap L$. A local leaf is called {\it univalent over} $V$, if it contains no singular 
 points and is bijectively projected onto $V$. 
 \end{definition}

  \begin{corollary} \label{cconva} 
 Let $a$ be a repelling fixed point of $f$, $\ha\in\caf^l$ be its fixed orbit. Let 
 $V=V(a)\subset\oc$ be a neighborhood of $a$, 
 $\{\hb^m\}_{m\in\nn}$ be a sequence of points in $\caf$ such that $\pi_0(\hb^m)=a$ and 
 the local leaves $L(\hb^m,V)$ are univalent over $V$ (see the above definition). 
 Then $\hf^m(\hb^m)\to\ha$, as $m\to+\infty$. 
\end{corollary}

\begin{proof} For simplicity we assume that $\hb^m\in\caf^l$ (the opposite case is 
treated analogously with minor modifications, taking into account Remark \ref{remcaf}). 
By Proposition \ref{conva}, for the  proof of the convergence $\hf^m(\hb^m)\to\ha$ 
it suffices to prove the statements of  the proposition for  $\hz=\ha$ and 
$\hz^m=\hf^m(\hb^m)$. Its first statement, which  says that 
$\pi_{-j}(\hf^m(\hb^m))\to a$, as $m\to\infty$, for all 
$j$, follows from construction: $\pi_{-j}(\hf^m(\hb^m))=\pi_0(\hf^{m-j}(\hb^m))=
a$, whenever $m\geq j$.   Let us prove its  second statement 
on univalence. Without loss of generality we consider that the local leaf $L(\ha,V)$ is 
univalent over $V$ (one can achieve this by shrinking $V$). 
By definition, $\hf^{-N}(\ha)=\ha$ for any $N\in\mathbb N$. 
Fix an arbitrary domain $W\subset\overline{\cc}$, $a\in W$,  such that the local leaf $L(\ha,W)$ is univalent 
over $W$ and a subdomain $U$, $\overline U\subset W$, $a\in U$. 
It suffices to show that for 
any $N\in\nn$ 
the  local leaf $L(\hf^{m-N}(\hb^m),U)$ is univalent over $U$, whenever $m$ is large 
enough. To do this, 
consider the inverse branches $f^{-s}$ that fix $a$: they are analytic on $W$ 
and tend to $a$ uniformly on $\overline U$, as $s\to+\infty$  
(by definition and the Shrinking Lemma). Hence, any of these branches 
maps $W$ diffeomorphically onto its image. 
Therefore, there exists a $k\in\mathbb N$, 
$k\geq N$,  such that $f^{N-m}(\overline U)\subset V$ for any $m\geq k$. Fix 
an arbitrary $m\geq k\geq N$.  Then  $f^{N-m}:U\to f^{N-m}(U)$ is a conformal diffeomorphism. The local leaf 
$L(\hb^m,f^{N-m}(U))$ is univalent over $f^{N-m}(U)$. Indeed, 
it is contained in the local leaf $L(\hb^m,V)$, which is univalent over $V$ by 
the conditions of the corollary. Together with the above statement, this 
implies that the local leaf $L(\hf^{m-N}(\hb^m),U)$ is univalent over 
$U$. This proves  the corollary.
\end{proof}

\begin{definition} \label{defass2} 
A leaf of $\caf$ is {\it associated to a repelling (or parabolic) periodic point}  
if it is contained in $\caf^l$ and coincides with a leaf of $\caf^n$ 
that is associated to the above point (see Definition \ref{defass1}). 
In this case we also say that the corresponding leaves of $\chf$ and $\chff$ 
are associated to this point. 
\end{definition}

\begin{proposition} \label{minim} \cite{klr, lm} The laminations $\caf$ and $\chf$ are minimal (i.e., each leaf is dense), if and only if the function $f$ does not have 
branch-exceptional repelling periodic orbits (see Definition \ref{crcom}). 
If $f$ has branch-exceptional repelling periodic orbits, then 
each of the above laminations has a finite number of isolated leaves 
(at most four; all of  them are associated to the latter periodic orbits) and becomes 
minimal after removing the isolated leaves. The  isolated leaves 
accumulate to their (minimal) complement, $\caf'$ or  $\chf'$:   
\begin{equation}\caf'=\caf\setminus(\text{the isolated leaves}); \ 
\chf'=\pi_{hyp}^{-1}(\caf')=\chf\setminus(\text{the isolated leaves}).
\label{noniso}
\end{equation}
One has $\chf'=\chf$, if and only if $f$ does not have branch-exceptional 
repelling periodic orbits. 
\end{proposition}

\subsection{Background material 3: horospheres; metric properties and basic cocycle} 
The horospheres in the hyperbolic 3- space with a marked point ``infinity''  
at the boundary 
(and in the leaves of the hyperbolic laminations) were defined in Subsection 1.1. We use their following well-known equivalent definition. Consider the 
projection $\pi:\hh\to L=\partial\hh\setminus\infty$ to the boundary plane along the geodesics 
issued from the infinity. In the 
model of half-space this is the Euclidean orthogonal projection to the boundary plane. It 
coincides with the natural projection $\htr=Aff(\cc)\slash S^1\to\cc=Aff(\cc)\slash\cc^*$, and its latter 
description equips the boundary with a natural complex affine structure: $L=\cc$. 
The boundary admits a Euclidean affine metric (uniquely defined 
up to multiplication by constant). 

Everywhere below whenever we consider a Riemannian metric on a surface, we treat 
it as a length element, not as a quadratic form. If we say ``two metrics 
are proportional'', then by definition, 
the proportionality coefficient is the ratio of the corresponding length elements. 

Consider a global section of the  projection $\pi:\hh\to L$: a surface 
in $\hh$ that is 1-to-1 projected to $L$. It 
carries two metrics: the restriction to it of  the hyperbolic metric of the 
ambient space $\hh$; the pullback of the Euclidean metric of $L$ under the projection. 

\begin{definition} \label{defheight} 
The above section is a {\it horosphere}, if its latter (Euclidean) metric 
is obtained from the former one (the restricted hyperbolic metric) by multiplication by a 
constant factor.  The {\it height} of a horosphere (with respect to the chosen Euclidean 
metric on $L$) is the logarithm of the latter constant factor. The {\it height of a given 
 point}  in the hyperbolic space is the height of the horosphere that contains this point.  \end{definition} 
 
 \begin{remark} \label{metholo} The height is a real-valued analytic function $\htr\to\rr$. 
 In the upper half-space model the horospheres are horizontal planes. Their  heights 
 (with respect to the standard Euclidean metric on $\partial\htr\setminus\infty$)  
  are equal to the logarithms of their Euclidean heights in the ambient  3- space. 
 The isometric liftings to $\htr$ of the affine mappings 
$z\mapsto\lambda z+b$ of the boundary $\cc=\partial\htr\setminus\infty$ transform 
the horospheres to the horospheres so that the height of the  image equals $\ln|\lambda|$ plus 
the height of the preimage. The horospheres foliate $\htr$. A quotient of 
$\htr$ by a discrete group of isometries fixing $\infty$ (e.g., a leaf of 
$\chf$ or $\chff$) carries the quotient foliation by horospheres. 
\end{remark}
 
 Now we discuss metric properties of the horospheres in the hyperbolic laminations. 
 Let $\caf$, $\chf$ be respectively the affine and the hyperbolic laminations 
 associated to a rational function $f$. 
 
\begin{proposition} \label{horolam} 
The horospheres form a lamination of any of the spaces $\chf$, 
$\chff$ by affine surfaces (orbifolds that may have conical singularities). In particular, 
the closure of a union of horospheres is a union of horospheres. 
\end{proposition}
 The proposition follows from the continuity of the hyperbolic metrics 
  (Remark \ref{afcont}). 
 
 Let $L\subset\caf$ be a leaf, $\hz\in L$ be a nonsingular point  such that the 
 restricted projection $\pi_0|_L$ has nonzero derivative at $\hz$. Fix a Hermitian metric on 
the tangent line to $\oc$ at $\pi_0(\hz)$. Its  projection pullback to the tangent line 
$T_{\hz}L$ extends (in unique way) up to a  Euclidean affine metric on the whole leaf $L$. 
Let $H$ be the corresponding leaf in $\chf$. We set  

\begin{equation}\beta_{\hz}:H\to\rr \ \text{the height with respect to the latter metric on } 
L, \text{ see Definition \ref{defheight}},\label{defhei}
\end{equation}
$$\a=(\hz,h)\in H \ \text{the point such that} \ \pi_{hyp}(\a)=\hz \ \text{and} \ 
\beta_{\hz}(\a)=h,$$
then we say that the point $\a$ is {\it situated over} $\hz$ at height $h$,
\begin{equation}S_{\hz,h}\subset H \ \text{the horosphere such that} \ 
\beta_{\hz}|_{S_{\hz,h}}\equiv h.\label{horoh}
\end{equation}

 \begin{proposition} \label{convh}  A sequence of points $(\hz^k,h_k)\in\chf$ converges to a point 
 $(\hz,h)\in\chf$, if and only if $\hz^k\to\hz$ in $\caf$ and $h_k\to h$. 
 \end{proposition}

The proposition follows from the continuity of the 
hyperbolic metrics of the leaves.  

When we extend the horospheres along loops in $\oc$, their heights may change. 
The monodromy of the heights is described by basic cocycle. Let us recall its definition. 

\begin{definition} \label{defbus}  Let $L\subset\caf$ be a leaf, 
$\hz,\hz'\in L$ be a pair of nonsingular points projected  to one and the same 
$z=\pi_0(\hz)=\pi_0(\hz')\in\oc$ so that the restricted projection $\pi_0|_L$ has nonzero 
derivative at both points $\hz$ and $\hz'$. Let $H=H(\hz)\subset\chf$ be the corresponding 
hyperbolic leaf. 
Fix  a Hermitian metric on $T_z\oc$. Let $\beta_{\hz},\beta_{\hz'}:H\to\rr$ be the corresponding heights defined 
 in (\ref{defhei}).  The {\it basic cocycle} is the difference 
$$\beta(\hz,\hz')=\beta_{\hz'}-\beta_{\hz}.$$
\end{definition}

\begin{remark} In the conditions of the above definition the basic cocycle is a 
well-defined constant and  
depends only on $\hz$ and $\hz'$ (it is independent on the choice of metric). One has 
$$\beta(\hz,\hz)=0, \ \b(\hz,\hz')=-\b(\hz',\hz).$$ 
Each horosphere $S_{\hz,h}\subset H(\hz)$ coincides with 
the horosphere $S_{\hz',h+\beta(\hz,\hz')}$. For any triple of nonsingular 
points $\hz,\hz',\hz''\in\caf$ lying in one and the same leaf $L$ and 
projected by $\pi_0|_L$ to one and the same point $z\in\oc$ with nonzero 
derivatives one has 
\begin{equation}\b(\hz',\hz'')=\b(\hz,\hz'')-\b(\hz,\hz') \ \text{(the cocycle 
identity)}.\label{cocycle}
\end{equation}The basic cocycle is $\hf$- invariant:
\begin{equation}\b(\hz,\hz')=\b(\hf^n(\hz),\hf^n(\hz')) \ \text{for any} \ n\in\mathbb N.\label{businv} 
\end{equation} 

\end{remark}

The next proposition is well-known and follows immediately from definition. 
\begin{proposition} Let $L\subset\caf$ be a leaf, $\hat c, \hat c'\in L$, 
$\pi_0(\hat c)=\pi_0(\hat c')=c$. Let $V\subset\oc$ be a neighborhood of $c$ such that the local leaves $L(\hat c,V),L(\hat c',V)\subset L$ are univalent over $V$ (see 
Definition \ref{univdefn}). Consider the conformal isomorphism 
\begin{equation}\psi_{\hat c,\hat c'}=(\pi_0|_{L(\hat c',V)})^{-1}\circ\pi_0|_{L(\hat c,V)}:
L(\hat c,V)\to L(\hat c',V).\label{psicc'}\end{equation}
Let us fix a Euclidean affine metric on the leaf $L$, which contains the above local 
leaves. Consider the module 
$|\psi'_{\hat c,\hat c'}|$ of derivative in the chosen Euclidean metric. 
Then for any $\hz\in L(\hat c, V)$, set 
$\hz'=\psi_{\hat c,\hat c'}(\hz)$, one has 
\begin{equation} \b(\hz,\hz')=-\ln|\psi'_{\hat c,\hat c'}(\hz)|.\label{busfun}\end{equation}
\end{proposition}

\begin{corollary} \label{remhar} Let $L$, $\hat c$, $\hat c'$, $V$ be as in the above proposition. 
For any $z\in V$ set 
$$\hz=\pi_0^{-1}(z)\cap L(\hat c,V), \ \hz'=\pi_0^{-1}(z)\cap L(\hat c',V). \ \text{The function}$$ 
\begin{equation}\b_{\hat c,\hat  c'}(z)=\b(\hz,\hz')\label{bushar}\end{equation}
 is harmonic on $V$, and hence, real-analytic. 
\end{corollary}

\subsection{Main results} 
First let us recall the following 
\begin{theorem} \label{lylist} (\cite{ly}, p.62) The affine lamination $\caf$ 
associated to a rational function $f$ 
admits a continuous family of Euclidean affine metrics on the non-isolated 
leaves, if and only if  $f$ is conformally-conjugated to a function from the 
list (\ref{list}). In the latter case there exists a unique 
(up to multiplication by constant) 
conformal Euclidean metric on $\oc$ with isolated singularities whose 
pullback to the non-isolated leaves under 
the projection $\pi_0:\caf'\to\oc$ yields the above Euclidean metric family on 
the non-isolated leaves.
\end{theorem}

\begin{corollary} \label{nodense} Let $f$ be a rational function from 
(\ref{list}). Then each horosphere in any non-isolated leaf 
of its quotient hyperbolic lamination $\chff$ is nowhere dense in $\chff$.
\end{corollary}

\begin{proof} Let $S$ be an arbitrary  horosphere in a non-isolated leaf of 
$\chf$. For the proof of the corollary it suffices to show that the images of 
$S$ under forward and backward iterations of $\hf$ are nowhere dense. Let  
$g$ denote the singular Euclidean metric on $\oc$ from the above theorem. 
We measure the heights of the horospheres with respect to the metric $g$ lifted 
to the leaves of $\caf'$ under the projection $\pi_0$. 
The heights of $S$ over all the points 
of the corresponding leaf of the affine lamination are all the same 
(by definition and Theorem \ref{lylist}). The mapping $f$ has a constant module 
of derivative in the metric $g$, since $\hf$ is leafwise affine. The heights 
of the iterated images of $S$ form  an arithmetic progression with step equal to 
the logarithm of the latter module of derivative. This 
progression  is  a discrete set of real numbers. Hence, the union of the  
 images of $S$ is nowhere dense. This proves the corollary. 
\end{proof}

\def\haf{H(\ha)\slash\hf}

Recall that $\chf'$ denotes the union of  the non-isolated leaves in $\chf$. 
It is $\hf$- invariant, and its quotient $\chf'\slash\hf$ is the  union of the non-isolated 
leaves in $\chf\slash\hf$. 

\begin{theorem} \label{trep}  
Let $f$ be a rational function that does not belong to the list 
(\ref{list}). Let $\chff$ ($\chf'\slash\hf$) be the corresponding quotient 
hyperbolic lamination 
(with deleted isolated leaves, see (\ref{noniso})). Let    
$H\slash\hf\subset\chf'\slash\hf$ be a non-isolated leaf associated to a repelling periodic 
point of $f$ (see Definition \ref{defass2}). 
Then each horosphere in $H\slash\hf$ is dense in $\chf'\slash\hf$. 
\end{theorem}

\begin{remark} \label{rtrep} The latter repelling periodic point is not branch-exceptional  
(Definition \ref{crcom}, Proposition \ref{minim} and Corollary \ref{coruniq}).
\end{remark}

Theorem \ref{trep} is the principal result of the paper. It is proved in Section 2. 
As it is shown below, it implies density of all the horospheres in $\chffp$ 
in the critically-nonrecurrent nonparabolic case 
and density of ``almost'' all the horospheres in the general 
critically-nonrecurrent case, provided that $f\notin\text{(\ref{list})}$. 
\def\wt#1{\widetilde#1}
\begin{theorem} \label{thyp} Let $f:\oc\to\oc$ be a critically-nonrecurrent 
rational function without parabolic periodic points (e.g., a hyperbolic one) that does 
not belong to the list (\ref{list}). Then each horosphere in  $\chff$ 
accumulates to $\chffp$.
\end{theorem}
\begin{theorem} \label{tnrec} Let $f$ be a critically-nonrecurrent rational 
function that does not belong to the list (\ref{list}). 
Let $H\subset\chf$ be  an arbitrary hyperbolic  leaf, 
$L=\pi_{hyp}(H)\subset\caf$ be the corresponding affine leaf. Let 
the projection $p(L)\subset \caf^n$ do not lie in a  leaf associated to a parabolic 
periodic point of $f$. Let $H\slash\hf\subset\chff$ be the corresponding leaf of the quotient lamination. 
Then each horosphere in $H\slash\hf$ accumulates to $\chffp$. 
\end{theorem}

Theorem \ref{thyp} follows immediately from Theorem \ref{tnrec}. Below we deduce 
Theorem \ref{tnrec}  from Theorem \ref{trep} and the following theorem, 
which will be proved in Section 3.  

\begin{theorem} \label{taccum} Let the conditions of Theorem \ref{tnrec} hold, but now 
$f$ is not necessarily excluded from the list (\ref{list}). 
Then each horosphere in $H\slash\hf$ accumulates to some horosphere in a leaf in 
$\chff$ associated to an appropriate repelling periodic point that is not 
branch-exceptional. 
\end{theorem} 

\begin{proof} {\bf of Theorem \ref{tnrec}.} Each horosphere in $H\slash\hf$ 
accumulates to some horosphere in 
a leaf in $\chff$ corresponding to some repelling periodic point that is not 
branch-exceptional (Theorem  \ref{taccum}). 
The latter horosphere accumulates to $\chffp$ (Theorem \ref{trep}). Hence, so 
does the former horosphere. This proves Theorems \ref{tnrec} and \ref{thyp}. 
\end{proof}

The following theorem proves the converse for the horospheres in the leaves 
associated to parabolic periodic points, without the critical nonrecurrence 
assumption.  

\begin{theorem} \label{taccum2} 
Let $f:\oc\to\oc$ be an arbitrary rational function with a parabolic 
periodic point $a$. Let $H_a\subset\chf$ be a leaf associated to it, 
$H_a\slash\hf\subset\chf\slash\hf$ be 
the corresponding leaf of the quotient hyperbolic lamination. Each  
horosphere in $H_a$ ($H_a\slash\hf$)  is closed in $\chf$ 
(respectively, $\chff$) and does not accumulate to itself. 
\end{theorem}

Theorem \ref{taccum2} will be proved in Section 4.

\def\fe{f_{\var}}
\def\chfe{\Cal H_{\fe}}
\def\be{B_{\var}}
\def\th{\tilde h}
\def\hfe{\hat f_{\var}}

The next theorem shows that  Theorems \ref{trep}-\ref{tnrec} 
are false when formulated for the 
nonfactorized lamination $\chf$. It deals with the quadratic polynomial family 
$$f_{\var}(x)=x^2+\var, \ \var<\frac14. \ \text{Set}$$
\begin{equation}a(\var)=\frac{1+\sqrt{1-4\var}}2: \ \text{this is 
the maximal real fixed point of} \ f_{\var}.
\label{aedef}\end{equation}

\begin{remark} \label{rtche} A polynomial $\fe$ belongs to the list (\ref{list}), if and only if 
either $\var=0$, or $\var=-2$. In the latter case $\fe$ is conformally conjugated 
to a Chebyshev polynomial. 
\end{remark}

\begin{proposition} \label{prex} The fixed point $a(\var)$ is repelling. 
It is branch-exceptional, if and only if $\var=-2$. 
\end{proposition}
\begin{proof} Recall that $\var<\frac14$. The fixed point equation $x^2-x+\var=0$ has two distinct solutions, one of them is $a(\var)$. 
Therefore, $\fe'(a(\var))\neq1$. One has $a(\var)>0$, by (\ref{aedef}). 
Thus, $a(\var)$ is not a parabolic fixed point: $1\neq f'(a(\var))=2a(\var)>0$. 
On the other hand, the polynomial $\fe$ has no fixed points 
on the right from $a(\var)$ by definition. It also has no critical points 
there, since $f'(x)=2x>2a(\var)>0$ for any $x>a(\var)$. 
This implies that the orbit of each point 
$x>a(\var)$ converges to infinity, since this is true whenever 
$x$ is large enough. Hence, $a(\var)$ is repelling. The quadratic 
mapping $\fe$ has unique 
critical point 0, and $\fe^{-1}(a(\var))=\{\pm a(\var)\}$. Hence, the fixed point $a(\var)$ is 
branch-exceptional, if and only if $-a(\var)$ is the critical value, i.e., $-a(\var)=\var$. 
The latter equation has the unique real solution $\var=-2$. 
\end{proof}

\def\cafe{\mathcal A_{\fe}}

\begin{theorem} \label{tnodense} Let $\var\in(-\infty,\frac14)\setminus
\{-2,0\}$, $\fe$, $a(\var)$ be as above, 
$\ha\in\cafe$ be the fixed orbit of $a(\var)$, $H(\ha)\subset\chfe$ be the 
corresponding leaf. Then each 
horosphere in $H(\ha)$ is nowhere dense in $\chfe$. More precisely, if 
$\pm\var>0$, then 
there exists a countable closed (and thus, nowhere dense) 
additive semigroup $\be\subset\rr_{\pm}$ such that for any $h\in\rr$ the 
accumulation set of the horosphere $S_{\ha,h}\subset H(\ha)$ is the horosphere union 
$$\cup_{\b\in\be}S_{\ha,h+\b}\subset H(\ha).$$\end{theorem}

{\bf Addendum.} {\it There exists an open set $W\subset\cc$ containing the 
interval union 
$(-\infty,0)\cup(0,\frac14)$ such that the statements of Theorem \ref{tnodense} hold 
true for all $\var\in W\setminus\{-2\}$, with $\pm\var>0$ replaced by 
$\pm\operatorname{Re}\var>0$.} 

\medskip

Theorem \ref{tnodense} and its addendum will be proved in Section 5.

\section{Density of the horospheres over repellers. Proof of Theorem \ref{trep}} 

\subsection{The plan of the proof of Theorem \ref{trep}} 
Let $a\in\oc$ be some repelling periodic point of $f$ that is not 
branch-exceptional (see Remark \ref{rtrep}). Let $\ha\in\caf$ be its  
periodic backward orbit, $L(\ha)$, $H(\ha)$ be the respectively the \
corresponding leaves of the laminations $\caf$ and $\chf$. Then the leaves 
$L(\ha)$, $H(\ha)$ are contained respectively in $\caf'$ and $\chf'$, and 
they are dense there (Proposition \ref{minim} and Corollary \ref{coruniq}). 
We fix a  horosphere $S\subset H(\ha)$, set  
 \begin{equation}\Cal S=\cup_{m\in\zz}\hf^m(S), \ \text{and show that the closure of} \ 
 \Cal S \ \text{in} \ \chf \ \text{contains} \  H(\ha).\label{denseinha}\end{equation}
 Then $\mathcal S$ is dense in $\chf'$, as is $H(\ha)$. This proves Theorem \ref{trep}. 
 
 It suffices to prove (\ref{denseinha}) with $S=S_{\ha,0}$. This implies the 
same statement for any other horosphere $S_{\ha,h}$, since for any $\hy\in L(\ha)$,  
$h\in\rr$ and $m\in\mathbb Z$ 
\begin{equation}  
(\text{the height of } \hf^m(S_{\ha,h}) \text{ over } \hy) \ = \ h+ 
(\text{the height of } \hf^m(S_{\ha,0}) \text{ over } \hy).\label{heq}
\end{equation}

 Without loss of generality everywhere below we assume that the point $a$ 
 is fixed: 
 $f(a)=a$.  One can achieve this by replacing $f$ by its  iterate: statement 
 (\ref{denseinha}) for an iterate, which we will prove, is stronger 
 than that for $f$. Then both leaves 
 $L(\ha)$ and $H(\ha)$ are fixed by $\hf$, which acts on $L(\ha)$ by  (complex) homothety 
 centered at $\ha$ with coefficient $f'(a)$. Set  
  \begin{equation}\Pi_a=\{\hy\in L(\ha)\setminus\ha\ | \ \pi_0(\hy)=a, \ (\pi_0|_{L(\ha)})'(\hy)\neq0\}.
  \label{univa}\end{equation}
 Each horosphere $S\subset H(\ha)$ is mapped by $\hf$ to a horosphere in the same  leaf $H(\ha)$ 
 so that
 \begin{equation} \hf^m(S_{\ha,0})=S_{\ha,m\ln|f'(a)|}, \ 
 \hf(S_{\hy,h})=S_{\hf(\hy),h+\ln|f'(a)|} \ \text{for any} \ m\in\zz, \ h\in\rr \ \text{and} \ \hy\in\Pi_a.\label{fsa}
 \end{equation}
  The monodromies of the horospheres  (when defined) along loops based at $a$ 
  add appropriate basic cocycles to the heights (see Definition \ref{defbus}) so that for any 
  $\hy\in\Pi_a$, $h\in\rr$, $m\in\zz$ 
  \begin{equation}
  S_{\ha,h}=S_{\hy,h+\b(\ha,\hy)}, \ \text{thus,}\ \hf^m(S_{\ha,0})=S_{\hy,h_{\hy,m}}, \ h_{\hy,m}=
  \b(\ha,\hy)+m\ln|f'(a)|.\label{monodr}\end{equation}
  The main part of the proof of Theorem \ref{trep} 
  is the next lemma, which implies that the above height values $h_{\hy,m}$ are dense in $\rr$. 
  Theorem \ref{trep} is then  deduced from it by elementary topological 
  arguments (using Corollary \ref{cconva}), which are presented at the end of the subsection. 
  
   \begin{proposition} \label{pianep} Let $f$ be a rational function, $a$ be some 
   its repelling  fixed point that is not branch-exceptional. 
   Then the set $\Pi_a$ from (\ref{univa}) is nonempty. 
   \end{proposition}
   \begin{proof} The 
   leaf $L(\ha)$ is contained in $\caf'$, is dense there, 
   and in particular, accumulates to itself (see the beginning of the  subsection).  
   Therefore, there exist a neighborhood $U=U(a)$ and  a sequence of points 
   $\ha^n\in L(\ha)$ converging to $\ha$ such that all the local leaves 
   $L(\ha,U)$, $L(\ha^n,U)$ are univalent and distinct (see Proposition \ref{conva}). 
   Then without loss of generality we consider that  
   $a^n_0=a$. By construction,  $\ha^n\neq\ha$ for infinite number of indices 
   $n$, and hence, $\ha^n\in\Pi_a$ (the above univalence statement). 
\end{proof}
   
 \begin{lemma} \label{vysoty} Let $f$ be a rational function that does not belong 
 to the list (\ref{list}), $a$ be some its repelling fixed point that is 
 not branch-exceptional. Let $\Pi_a$ be as in (\ref{univa}). 
 The set 
 \begin{equation} \Cal Bf=\{ \b(\ha,\hy)+m\ln|f'(a)| \ | \ \hy\in\Pi_a, \ m\in\zz\}\label{cab}
 \end{equation} 
 is dense in $\rr$.
  \end{lemma}
 Everywhere below for any $z\in\oc$ (with a chosen local holomorphic 
 chart in its neighborhood, 
 the latter being equipped  with the standard Euclidean metric) and $\delta>0$ 
 we set  
 $$D_{\delta}(z)=\text{ the } \delta- \text{ disk centered at } z, 
 \ D_{\delta}=D_{\delta}(0).$$
 \def\dda{D_{\delta}(a)}
  
 The proof of  Lemma \ref{vysoty} modulo technical details is given below. The details 
 of the proof take the most part of the section. In its proof 
  we use the following properties of the points from $\Pi_a$ and basic cocycles. 

 \begin{proposition} \label{psumbus} Let $f$ be a rational function, $a\in\oc$ 
 be some its repelling fixed point that is not branch-exceptional. Let 
 $\Pi_a$ be as in (\ref{univa}), $\hy,\hc\in\Pi_a$. Let $\delta>0$ be such that the local 
 leaves $L(\ha,\dda)$, $L(\hy,\dda)$, $L(\hc,\dda)$ are  univalent over $\dda$, and 
 moreover, the inverse branch $f^{-1}$ that fixes $a$ extends up to a univalent holomorphic 
 function $\dda\to\dda$ 
 (whose orbits in $\dda$ thus converge to $a$). Let $j\in\mathbb N$ be such that $y_{-k}\in\dda$ 
 for any $k\geq j$ (see Proposition \ref{repell}). Let 
 $$\hat\zeta\in L(\hc,\dda), \ \pi_0(\hat\zeta)=y_{-j},\ \hw=\hf^j(\hat\zeta), \ \b_{\ha,\hc} \ 
 \text{be the function from (\ref{bushar}). Then}$$
 \begin{equation}\hw\in\Pi_a \text{ and } 
 \b(\ha,\hw)=\b(\ha,\hy)+\b_{\ha,\hc}(y_{-j}).\label{sumbus}
 \end{equation}
 \end{proposition}
  
  \begin{proof} One has $\hw\in\Pi_a$. Indeed, 
  $\pi_0(\hw)=f^j(y_{-j})=a$. No $w_{-k}$ is a critical point of $f$. 
 For $k> j$ this 
 follows from definition and the univalence of the local leaf 
 $L(\hc,D_{\delta}(a))$. 
 For $k\leq j$ one has $w_{-k}=y_{-k}$. The latters are not critical 
 points of $f$, since $\hy\in\Pi_a$. Hence, $\hw\in\Pi_a$.  
 
 By (\ref{cocycle}), 
 \begin{equation}\b(\ha,\hw)=\b(\ha,\hy)+\b(\hy,\hw).\label{cosum}\end{equation}
 We show that 
 \begin{equation}\b(\hy,\hw)=\b_{\ha,\hc}(y_{-j}).\label{baaa}\end{equation}
 This together with (\ref{cosum}) implies (\ref{sumbus}). One has 
  $$\b(\hy,\hw)=\b(\hf^{-j}(\hy),\hf^{-j}(\hw))=\b(\hat t,\hat\zeta), \ 
  \hat t=\hf^{-j}(\hy),$$  
 by (\ref{businv}) and since $\pi_0(\hat t)=\pi_0(\hat\zeta)=y_{-j}$. 
 The point $\hat t$ lies in the local leaf $L(\ha,\dda)$. This follows 
 from the choice of $j$ ($t_{-s}=y_{-j-s}\in\dda$ for any $s\geq0$) 
 and the invariance of $D_{\delta}(a)$ under the inverse branch $f^{-1}$ 
 fixing $a$ (the condition of the proposition). Thus, 
 $$\hat t\in L(\ha,\dda), \ \hat\zeta\in L(\hc,\dda), \ \pi_0(\hat t)=\pi_0(\hat\zeta)=y_{-j}. \ 
 \text{Therefore,}$$
 $$\b(\hy,\hw)=\b(\hat t,\hat\zeta)=\b_{\ha,\hc}(y_{-j}),$$
by (\ref{bushar}). This  proves (\ref{baaa}) and (\ref{sumbus}).
 \end{proof} 
 
 \begin{corollary} \label{semigr} Let $f$, $a$, $\Pi_a$ be as in Proposition \ref{psumbus}. 
  The closure 
 \begin{equation}\Cal B=\overline{\{ \b(\ha,\hy)\ | \ \hy\in\Pi_a\}} \label{calb}\end{equation}
  is an additive semigroup in $\rr$. 
   \end{corollary}
   
   \begin{proof}
We have to show that for any given $\hy,\hc\in\Pi_a$ and $\var>0$ there exists a 
$\hw\in\Pi_a$ such that $|\b(\ha,\hw)-(\b(\ha,\hy)+\b(\ha,\hc))|\leq\var$. 
Let $y_{-j}$ and $\hw$ be as in the above proposition. Then $\hw$ is a one 
we are looking for, whenever $j$ is large enough. Indeed, the previous difference equals 
$\b_{\ha,\hc}(y_{-j})-\b_{\ha,\hc}(a)$, by (\ref{sumbus}) and since 
$\b(\ha,\hc)=\b_{\ha,\hc}(a)$ by definition. The latter difference tends to 0, as 
$j\to+\infty$, since $y_{-j}\to a$. This  proves the corollary.
\end{proof}

   We use the following elementary property of additive semigroups. 
\begin{proposition} \label{pairs} 
Let $\Cal B\subset\rr$ be an additive semigroup such that for any $\var>0$ 
 it contains a pair of at most $\var$- close  distinct elements. Then for any $M\in\rr\setminus0$ 
 the semigroup $\Cal B_M=\Cal B+\zz M$ is dense in $\rr$.
 \end{proposition} 

\begin{proof} Fix a $\var>0$ and a pair $A,B\in\Cal B$, $A\neq B$, 
 $\var'=|A-B|\leq\var$. The elements $A$ and $B$ generate a semigroup that contains 
 arithmetic progressions with step $\var'$ and arbitrarily large lengths: for any $m\in\mathbb N$ 
 the set $\{ sA+(m-s)B \ | \ 0\leq s\leq m\}$ is such a progression. Fix a $m$ 
 such that the ends $mA$, $mB$ of the latter progression bound an interval of 
 length at least $M$. The images of the progression by translations by $nM$, 
 $n\in\mathbb Z$, form a $\var'$- net on $\rr$. Therefore, the set 
 $\Cal B_M$ from the proposition contains 
 a $\var'$- net on $\rr$ with arbitrarily small $\var'$. Hence, it is dense. The proposition 
 is proved.
 \end{proof}

By definition, one has 
\begin{equation}\Cal B+\zz\ln|f'(a)|\subset\overline{\Cal Bf}.\label{newincl}
\end{equation}
We show that the semigroup  $\Cal B$ contains 
distinct elements arbitrarily close to each other. Then applying Proposition \ref{pairs} to 
$M=\ln|f'(a)|$ together with (\ref{newincl}) implies Lemma \ref{vysoty}. 

As it is shown below, the above statement on $\Cal B$ is implied by (\ref{sumbus}) and  the following 

\begin{lemma} \label{b12} {\bf (Main Technical Lemma)}  
 Let $f$ be a rational function that does not belong to the  list (\ref{list}). 
 Let  $a\in\oc$ be some its repelling fixed point that is not branch-exceptional, 
$\ha\in\caf$ be its fixed orbit, $\Pi_a$ be the set from (\ref{univa}).  
There exists a pair of points $\hy,\hc\in\Pi_a$ such that for any $N\in\mathbb N$ 
\begin{equation} \b_{\ha,\hc}|_{\{y_{-j} \ | \ j\geq N\}}\not\equiv const.\label{nonconst1}\end{equation}
\end{lemma}
  
 The proof of Lemma \ref{b12} (given in the rest of the section) uses essentially the analyticity of basic cocycle. 
 
 \begin{proof} {\bf of Lemma \ref{vysoty} modulo Lemma \ref{b12}.} It suffices to show that 
the semigroup $\Cal B$ contains 
 pairs of arbitrarily close distinct elements (see the above discussion).  
 Let $\hy,\hc\in\Pi_a$ be as in Lemma \ref{b12}. 
 The value $B'=\b(\ha,\hy)+\b(\ha,\hc)$ is contained in $\Cal B$ (Corollary \ref{semigr}). 
 On the  other hand, for any $j$ and $\hw$ as in (\ref{sumbus}) the value 
 $\b(\ha,\hw)$ also belongs to $\Cal B$. It differs from $B'$ by 
 $\b_{\ha,\hc}(y_{-j})-\b_{\ha,\hc}(a)$, which tends to 0, as $j\to+\infty$ 
 (see the proof of Corollary \ref{semigr}). The latter difference is nonzero for an 
 infinite number of 
 values of $j$ (Lemma \ref{b12}). Thus, the previously constructed elements $B'$ and 
 $\b(\ha,\hw)$  of the semigroup  $\Cal B$  can be made distinct 
 and arbitrarily close to each other. This proves Lemma \ref{vysoty}.
  \end{proof}

\begin{proof} {\bf of Theorem \ref{trep} modulo Lemma \ref{b12}.} To prove the density of the set $\Cal S$, see 
(\ref{denseinha}),  with $S=S_{\ha,0}$, we use Lemma \ref{vysoty} and 
Proposition \ref{horolam}.  
 For any $h\in\mathcal Bf$ we construct  a sequence $\{\ha^m\}|_{m\in\nn}
 \subset\Pi_a$ such that $(\ha^m,h)\in\mathcal S$ for all $m$ and 
 $(\ha^m, h)\to(\ha, h)$. This together with the density of 
 $\mathcal Bf$ (Lemma \ref{vysoty}) implies that $\mathcal S$ accumulates to all the 
  points  
 $(\ha,r)$, $r\in\rr$, and hence, to the horospheres through these points 
  (Proposition \ref{horolam}).  The latter horospheres saturate 
 the whole leaf $H(\ha)$, thus, $H(\ha)\subset\overline{\mathcal S}$. This together with 
 the discussion from the beginning of the section  implies Theorem \ref{trep}. 

Fix a $h\in\mathcal Bf$. By definition, there exist a 
  $\hy\in\Pi_a$ and a $l\in\zz$ such that $h+l\ln|f'(a)|=\beta(\ha,\hy)$, i.e., 
  $(\hy,h+l\ln|f'(a)|)\in S_{\ha,0}$, see (\ref{cab}). Then for any $m\in\nn$ 
  $$q_m=(\hy, h+(l-m)\ln|f'(a)|)\in S_{\ha,-m\ln|f'(a)|}\subset\Cal S,$$ 
 by (\ref{heq}), where $m$ is replaced by 0 and  $h$ is replaced by $-m\ln|f'(a)|$. 
 For any $m\geq l$ set 
  $$  
  \ha^m=\hf^{m-l}(\hy), \alpha_m=(\ha^m,h)=\hf^{m-l}(q_m).$$
  The sequence $\{\ha^m\}_{m\geq l}$ is a one we are looking for. 
  Indeed, $\alpha_m\in\Cal S$, since $\Cal S$ is $\hf$- invariant and 
$q_m\in\Cal S$. One has $a_0^m=f^{m-l}(y_0)=f^{m-l}(a)=a$, whenever $m\geq l$. 
Set $\hb^m=\hy$. Let $V$ be a neighborhood of $a$ such that the local leaf 
$L(\hy,V)$ is univalent over $V$.  The sequence $\ha^m$ converges to $\ha$ in $\caf$, by Corollary \ref{cconva} applied to these $\hb^m$ and $V$. 
This together with Proposition \ref{convh} proves 
 the convergence $\alpha_m\to(\ha,h)$ and finishes the proof of Theorem \ref{trep}. 
 \end{proof}

\subsection{Proof of the Main  Lemma \ref{b12} modulo technical details} 
Let $f$ be a rational function. Let $a$ be some its repelling fixed point 
that is not branch-exceptional, $\ha\in\caf$ 
be  its fixed orbit, $\Pi_a$ be the set from (\ref{univa}). 
Fix a small  neighborhood 
$U$ of $a$ where $f$ is univalent and such that $f(U)\supset U$. 
The branch of $f^{-1}$ fixing $a$ is holomorphic on $U$, 
sends $U$ to itself and its iterations 
converge to $a$ uniformly on compact subsets of $U$. 
(In particular, the local leaf $L(\ha,U)$ is univalent over $U$.) Therefore, the linearizing 
coordinate of $f$ at $a$ extends up to a conformal  coordinate on $U$. In addition, 
we will assume for a technical reason 
that $U$ is convex (e.g., a disk) in the linearizing coordinate. 

\begin{definition} Let $f$, $\ha$, $U$ be as above. Let 
$\hz=(z_0,z_{-1},\dots)\in L(\ha)$, $N>0$ be such that 
$z_{-j}\in U$ for any $j\geq N$. (The number $N$ is not necessarily the minimal one 
with this property.) The backward orbit  
$z_{-N}, z_{-N-1},\dots$ is called a {\it tail of} $\hz$. 
If $N$ is the minimal number as above, then the tail is called {\it complete}. 
(Equivalently, a (complete) tail is a (maximal) 
backward orbit representing a point of $L(\ha,U)$.) 
\end{definition}

\begin{lemma} \label{enonconst} Let $f$ be a rational function that does not belong to 
the list (\ref{list}), $a$, $\Pi_a$ be as above. There exists a point 
$\hc\in\Pi_a$ such that
 \begin{equation}\b_{\ha,\hc}\not\equiv const \ \text{in a neighborhood of} \ a.\label{bneconst}
 \end{equation} 
 \end{lemma}
This lemma is proved in the next subsection. 
Fix a $\hc\in\Pi_a$ satisfying (\ref{bneconst}). We prove the statement of 
Lemma \ref{b12} for this $\hc$: there exists a $\hy\in\Pi_a$ such that 
$\beta_{\ha,\hc}\not\equiv const$ along any tail of $\hy$. 
Without loss of generality we consider that the 
local leaf $L(\hc,U)$ is univalent over $U$. One can achieve this by shrinking $U$. 
 Recall that the function $\b_{\ha,\hc}$ is real-analytic on $U$. We consider the 
 auxiliary analytic subset $A\subset U$: 
 $$A= \text{ the intersection of all the analytic subsets in } U \text{ that contain 
 some tail of each} \ \hy\in\Pi_a.$$
\begin{proposition} \label{pau} Let $f$, $a$, $U$,  be as at the beginning of the 
subsection, $A$ be as above. The set 
 $A$ is either the whole $U$, or a line interval through $a$ 
 in the linearizing chart of $f$ on $U$.
 \end{proposition}
  
  \begin{remark} If $A$ is a line interval, then $f'(a)\in\rr$. This follows from 
  definition and the $f$- invariance of the germ of $A$ at $a$ (see Claim 2 in 2.4  
  for a stronger  invariance statement). 
  \end{remark}
 Recall that the linearizing coordinate on $U$ lifts to the local leaf $L(\ha,U)$ 
 and extends up to a global affine coordinate on $L(\ha)$ (Example \ref{exrep}). 
 We consider the auxiliary conformal mapping  
\begin{equation}
\psi_{\hc}:U\to L(\hc,U)\subset L(\ha)=\cc, \ \psi_{\hc}=(\pi_0|_{L(\hc,U)})^{-1}: 
\text{ one has } \beta_{\ha,\hc}\equiv-\ln|\psi_{\hc}'|,\label{defpsi}
\end{equation}
the derivative is taken in the above coordinates. This follows from 
 (\ref{busfun}). 

    \begin{proposition}\label{afia}  In the conditions of Proposition 
    \ref{pau} let $A$ be a line interval. 
   Let $\hat A\subset L(\ha)$ denote the line whose intersection with the local leaf 
   $L(\ha,U)$ is projected to $A$. For any $\hc\in\Pi_a$ such that the local 
   leaf $L(\hc,U)$ is univalent over $U$ the mapping 
 $\psi_{\hc}$  from (\ref{defpsi}) sends $A$  to  $\hat A$. 
 \end{proposition}
 Proositions \ref{pau} and \ref{afia} will be proved in 2.4. 
 
\begin{proof} {\bf of Lemma \ref{b12} modulo Propositions \ref{pau} and \ref{afia}.} 
Let $f$ do not belong to the list (\ref{list}). Let 
$a$, $\ha$, $U$, $A$, $\hc$ be as above. We prove Lemma \ref{b12} by contradiction. 
Suppose the contrary: $\beta_{\ha,\hc}\equiv const$ along some tail of each 
$\hy\in\Pi_a$. 
Then $\beta_{\ha,\hc}|_A\equiv\beta_{\ha,\hc}(a)=\beta(\ha,\hc)$, by definition, 
analyticity and since all the tails have the common limit $a$. If $A=U$, then the 
latter identity contradicts (\ref{bneconst}). Thus, $A$ is a line  interval (Proposition 
\ref{pau}) and $\beta_{\ha,\hc}|_A\equiv const$.  Therefore, the conformal mapping 
$\psi_{\hc}$ from (\ref{defpsi}) sends the line $A$ to the line $\hat A$ 
(Proposition \ref{afia}), and its 
derivative in the linearizing chart  has constant module along $A$. The next proposition 
shows that $\psi_{\hc}$ is an affine mapping. 

  \begin{proposition} \label{phiaf} 
 Let $\psi$ be a conformal mapping of one domain of  $\cc$ onto another one. 
 Let $\psi$ map 
 a line interval $A$ to a line and the module of its derivative be constant along $A$. Then 
 $\psi$ is an affine mapping.
\end{proposition}

\begin{proof} There exists a complex affine mapping that coincides with $\psi$ on $A$  (the constance of the module of the derivative on $A$ and the linearity of $A$ and its image). Then it 
coincides with $\psi$ everywhere (the uniqueness of analytic extension). This proves 
the proposition.
\end{proof}

Thus, the mapping $\psi_{\hc}$ is affine. Therefore, its derivative is constant on $U$, 
and hence, $\beta_{\ha,\hc}\equiv const$ on $U$, by (\ref{defpsi}), - 
a contradiction to (\ref{bneconst}).  Lemma \ref{b12} is proved.
\end{proof} 

\begin{remark} The above arguments show that if $\beta_{\ha,\hc}|_A\equiv const$, 
then $\beta_{\ha,\hc}\equiv const$. 
Thus, if the former holds true for all $\hc\in\Pi_a$, then 
the function $f$ belongs to the list (\ref{list}) (Lemma \ref{enonconst}). 
Generically, $A=U$.
There exist  functions $f$ that do not belong to the list (\ref{list}) but for which $A$ 
is a line interval, see the next example.
 \end{remark}

\begin{example} 
  Consider a  quadratic polynomial $f(x)=cx(1-x)$, $c>4$. It has a repelling 
 fixed point $a=0$. It is well-known that there exists a Cantor set in $[0,1]$ 
 containing 0 and completely invariant (i.e., coinciding with its preimage) under the 
 complex dynamics of $f$. The corresponding set $A$ is an interval of the 
 real line, since each complex preimage of 0 under arbitrary iterate of  $f$  is real: 
 it is contained in the above Cantor set. On the other hand, the polynomial $f$ does not belong to the list (\ref{list}). 
   \end{example}

\subsection{Nonconstance of basic cocycles. Proof of  Lemma \ref{enonconst}}
Fix an affine Euclidean metric $g$ on $L(\ha)$.  
We prove the lemma by contradiction. Suppose the contrary: 
$\b_{\ha,\hc}\equiv const$ in a neighborhood of $a$ for all $\hc\in\Pi_a$. We show 
(the next proposition) that the latter constant is zero. Afterwards we show 
(Proposition \ref{zerozero})
that the metric $g$ is projected to some well-defined (singular) metric on $\oc$ 
and extends up to a continuous family of affine Euclidean metrics on the 
leaves of $\caf'$. Hence, $f$ belongs to the list (\ref{list}), by 
Theorem \ref{lylist}, - a contradiction to the conditions of 
Lemma \ref{enonconst}. 

\begin{proposition} \label{zeroconst} Let $f$, $a$, $\ha$, 
$\Pi_a$ be as at the beginning of the previous subsection. Let $\hc\in\Pi_a$ 
be such that $\b_{\ha,\hc}\equiv const$ in a neighborhood of $a$. Then 
$\b_{\ha,\hc}\equiv0$. 
\end{proposition}

\begin{proof} Let $U$ be a neighborhood of $a$ such that the local leaves 
$L(\ha,U)$ and $L(\hc,U)$ are  univalent over $U$. 
Consider  the mapping 
$$\psi:L(\ha,U)\to L(\hc,U), \ \psi=(\pi_0|_{L(\hc,U)})^{-1}\circ\pi_0.$$
 The module of the derivative of 
$\psi$ in the Euclidean metric $g$ on $L(\ha)$ 
is constant, since its logarithm equals $-\b_{\ha,\hc}\equiv const$ 
by (\ref{busfun}). Therefore, $\psi$ extends up to an affine 
automorphism $\psi:L(\ha)\to L(\ha)$. If the module of its derivative equals 1, then 
$\b_{\ha,\hc}\equiv0$, and the statement of the 
proposition follows immediately. Otherwise $\psi$ has a fixed point (denote it $q$) 
that is either 
an attractor, or a repeller. The function $\pi_0:L(\ha)\to\oc$ is $\psi$- invariant by 
construction and analyticity. Therefore, 
it is constant in a neighborhood of $q$ (analyticity). Hence, $\pi_0|_{L(\ha)}\equiv const$, 
- a contradiction. The proposition is proved. 
\end{proof}

\begin{proposition} \label{zerozero} 
Let $f$, $a$, $\ha$, $\Pi_a$ be as at the beginning of the previous subsection. Let  
\begin{equation}\b_{\ha,\hc}\equiv0\ \text{for any} \ \hc\in\Pi_a. \label{eqmetr}\end{equation}
Consider arbitrary triple $(W,W_1,W_2)$, where $W\subset\oc$ is a simply \
connected domain  intersecting the Julia set of $f$, $W_1,W_2\subset L(\ha)$ 
are local leaves 1-to-1 projected to $W$ by $\pi_0$. Let $g$ be an affine 
Euclidean metric on $L(\ha)$. 
Then the pushforwards $(\pi_0)_*(g|_{W_i})$ of the metrics $g|_{W_i}$, 
$i=1,2$, to $W$ coincide. 
\end{proposition}

In the proof of the proposition we use the following

\begin{proposition} \label{brex} Let $f:\oc\to\oc$ be a rational function, $a\in\oc$ be 
some its 
repelling fixed point that is not branch-exceptional. Let $W\subset\oc$ be an arbitrary 
simply connected domain intersecting the Julia set. There exist a $n\in\mathbb N$ and 
a $z\in W$ such that $f^n(z)=a$ and $(f^n)'(z)\neq0$.
\end{proposition}

\begin{proof} Fix a nonperiodic backward orbit $\hat q=(q_0=a,q_{-1},\dots)$ 
avoiding critical points. Then $q_{-j}$ is not a postcritical point, whenever $j$ is large 
enough. Indeed, otherwise the nonperiodic (and hence, {\it infinite}) 
{\it backward} orbit $\hat q$ 
is contained in the finite union of {\it forward} postcritical orbits, which is impossible. 
Fix three distinct points $q_{-j_1}$, $q_{-j_2}$, 
$q_{-j_3}$, $j_1<j_2<j_3$. The union 
$\mathcal U=\cup_{m\in\mathbb N\cup 0}f^m(W)$ contains some 
of them. Otherwise, the family $f^m:W\to\oc$ of meromorphic functions avoids 
three distinct values and  is normal by Montel's theorem (\cite{lyu}, p.52). 
Hence, $W\cap J=\emptyset$, - a contradiction. This implies that 
there exist a $m\in\mathbb N$ and a $z\in W$ such that $f^{m}(z)$ coincides with 
some of the three above values. Hence, $f^n(z)=a$ for $n=m+j_3$. 
By construction, $(f^n)'(z)\neq0$, since no $q_{-j}$ is a critical point. 
This proves the proposition.
\end{proof}

\begin{proof} {\bf of Proposition \ref{zerozero}.} Let $z\in W$ and $n\in\mathbb N$ 
be as in Proposition \ref{brex},
$$\hz^i=(\pi_0|_{W_i})^{-1}(z)\in W_i, \ \hy^i=\hf^n(\hz^i), \ i=1,2.$$
Then  $\hy^i\in\Pi_a$. Indeed, $\hy^i\in L(\ha)$, since $\hz^i=\hf^{-n}(\hy^i)\in 
W_i\subset L(\ha)$. One has $\pi_0(\hy^i)=f^n(z)=a$. The germ of the 
projection $\pi_0:(L(\ha),\hy^i)\to(\oc,a)$ has nonzero derivative at $\hy^i$.  
Indeed, this germ is equal to the composition $f^n\circ\pi_0\circ\hf^{-n}$ of the 
following germs: 
$$\hf^{-n}:(L(\ha),\hy^i)\to(L(\ha),\hz^i), \ \ \pi_0:(L(\ha),\hz^i)
\to(\oc,z), \ \ f^n:(\oc,z)\to(\oc,a).$$
Each of them has nonzero derivative at the 
corresponding base point, since $\hf^{-n}$ is affine, each $\hz^i$ is contained in 
a univalent leaf $W_i$ and $(f^n)'(z)\neq0$. Thus, the derivative of the 
composition at $\hy^i$ is also nonzero. 
The statement of Proposition \ref{zerozero} is equivalent to the identity 
$\beta_{\hz^1,\hz^2}\equiv0$, by the definition of  basic cocycle. One has 
$\beta_{\hz^1,\hz^2}=\beta_{\hy^1,\hy^2}=\beta_{\ha,\hy^2}-\beta_{\ha,\hy^1}\equiv0$,  
by the invariance of basic cocycle, the cocycle identity and (\ref{eqmetr}). 
This proves Proposition \ref{zerozero}.
\end{proof} 

\begin{proof} {\bf of Lemma \ref{enonconst}.} Recall that we assume that 
$\beta_{\ha,\hc}\equiv const$, hence, $\beta_{\ha,\hc}\equiv0$ for all 
$\hc\in\Pi_a$ (Proposition \ref{zeroconst}). 
Let  $\hz,\hw\in L(\ha)$ be arbitrary two points 
with $z_0=w_0$ such that $(\pi_0|_{L(\ha)})'(\hz),
(\pi_0|_{L(\ha)})'(\hw)\neq0$. Consider the restrictions to $T_{\hz}L(\ha)$ 
and $T_{\hw}L(\ha)$ of the affine metric $g$ of the leaf $L(\ha)$. 
The latter are projected by $\pi_0$ to the same metric on $T_{z_0}\oc$. 
This follows by choosing a small neighborhood $W$ of $z_0$ such that the 
local  leaves $L(\hz,W)$ and $L(\hw,W)$ are univalent and applying
Proposition \ref{zerozero}. Those pairs of points $\hz$, $\hw$ 
form an open and dense 
subset in the space of pairs of points in $L(\ha)$ with coinciding projections. 
Hence, the metric $g$ of the whole leaf $L(\ha)$ is 
projected by $\pi_0$ to a well-defined 
(maybe singular) conformal metric on $\oc$. The only possible singularities 
of the projected metric are the so-called {\it persistently-critical values}: 
those $z\in\oc$, whose preimages under the projection $\pi_0|_{L(\ha)}$ 
either are empty or contain only critical points of the projection. 
 (There can be at most 4 persistently-critical values, by Ahlfors' 
 five island theorem, see \cite{ah5}, theorem VI.8.) 
 Lifting the projected metric to other leaves 
 extends $g$ up to a continuous family 
 of conformal metrics on all the leaves of $\caf'$ (that a priori 
 may have isolated singularities on some leaves).  The lifted metrics are 
 affine on all the leaves, as is $g$, by the density of the leaf $L(\ha)$ 
 (Proposition \ref{minim}) and the continuity of the affine structure family. 
 Therefore, $f$ belongs to the list (\ref{list}), - a contradiction.  
This proves Lemma \ref{enonconst}. 
\end{proof}

\subsection{The common level of basic cocycles. Proof of Propositions 
\ref{pau} and \ref{afia}}
Everywhere below we suppose that $A\neq U$. This implies that there exists 
a real-analytic curve containing a tail of each $\hy\in\Pi_a$. First we show that $A$ is a 
finite union of line intervals (Claim 1). Afterwards we prove Propositions \ref{pau} and 
\ref{afia}, using the invariance of $A$ under appropriate inverse branches 
(Claims 2 and 3). 

\medskip
\def\hq{\hat q}

{\bf Claim 1.} {\it The set $A$ is a finite union of line intervals intersecting at $a$ 
with ends in $\partial U$. It contains a complete tail of each $\hy\in\Pi_a$.}

\begin{proof} Any tail is an orbit under the 
iterations of the multiplication by $(f'(a))^{-1}$ in the linearizing chart. 
Therefore,  
$\arg f'(a)=\pi\frac pq\in\mathbb Q$: otherwise no tail can be 
contained in a real-analytic curve, and $A=U$, - a contradiction. 
Therefore, for any $\hy\in\Pi_a$ there exists a union $\Lambda_{\hy}$ of $q$ lines 
through $a$ such that each tail of $\hy$ is contained there and intersects 
each line at a sequence of points converging to $a$. The intersection 
$\Lambda_{\hy}\cap U$ is a union of $q$ intervals through $a$ with ends in 
$\partial U$, since $U$ is convex by assumption. For any $\hy\in\Pi_a$ the 
latter intersection lies in $A$ 
and contains the complete tail of $\hy$, by the two previous statements and 
analyticity. Thus, $A$ is a union of intervals through $a$. This union 
is finite by analyticity. Claim 1 is proved.  
 \end{proof}
 
 {\bf Claim 2.} {\it The set $A$ is invariant under the branch 
 $f^{-1}:U\to U$ fixing $a$.}
 
 \begin{proof} The $f^{-1}$- invariance of $A$ follows from definition and the 
 invariance of $U$.
 \end{proof}
 
 {\bf Claim 3.} {\it Let $n\in\mathbb N$, $f^{-n}$ be an arbitrary germ of 
  holomorphic inverse branch 
at $a$ such that $f^{-n}(a)\in U$. Then $f^{-n}$ sends the germ of $A$ at $a$ 
 to $A$.}
 
 \begin{proof} Fix a neighborhood $V\subset U$  of $a$ such that $f^{-n}$ is 
 holomorphic on $V$ and $f^{-n}(V)\subset U$. This defines the holomorphic branches 
 \begin{equation}f^{-j}=f^{n-j}\circ f^{-n}:V\to\oc \text{ for all } j\in\mathbb N:
 \label{invbr}\end{equation}
 for $j\leq n$ the function $f^{n-j}$ is well-defined on $\oc$; for $j>n$ we 
 take $f^{n-j}:U\to U$ to be the holomorphic branch fixing $a$. Then 
 $f^{-j}(V)\subset U$ for all $j\geq n$ and  
 $f^{-j}\to a$ uniformly on compact subsets in $V$, as $j\to+\infty$. 
 Therefore, $\hq=(a,f^{-1}(a),f^{-2}(a)\dots)\in\Pi_a$ and 
 $f^{-n}(a)$ is contained in the complete tail of $\hq$.  
Thus,  $f^{-n}(a)\in A$ (Claim 1). 
 
 Fix an arbitrary $\hy\in\Pi_a$. For any $l$ large enough $y_{-l}\in V$. 
 For these $l$ the values $f^{-j}(y_{-l})$ with $j\in\mathbb N$ are defined 
 by (\ref{invbr}). Those of them with  $j\geq n$ form a tail of the backward orbit 
 $$(y_0=a, y_{-1}, \dots, y_{-l}, f^{-1}(y_{-l}),f^{-2}(y_{-l}),\dots)\in\Pi_a.$$
 Therefore, $f^{-n}(y_{-l})\in A$ (Claim 1). This holds true for arbitrary 
 $\hy\in\Pi_a$ and any $l$ large enough (dependently on $\hy$). Thus, 
 $f^{-n}$ maps some tail of arbitrary $\hy\in\Pi_a$ to $A$. Hence, $f^{-n}$ sends 
 the germ of $A$ at $a$ to $A$, by definition and analyticity. This proves Claim 3.
 \end{proof} 
 
 \begin{proof} {\bf of Proposition \ref{pau}.} Suppose the contrary: $A\neq U$ and $A$  
 is not a single 
 interval. Hence, it is a finite union of line intervals through $a$ (Claim 1) and 
  contains at least two distinct line intervals $I_1$, $I_2$ 
   intersecting transversely at $a$. Fix an 
  arbitrary $\hy\in\Pi_a$. Let $(y_{-n}, y_{-n-1},\dots)$ be some its tail in $U$, 
  $f^{-n}$ be the corresponding inverse branch germ at $a$: $f^{-n}(a)=y_{-n}$. Then 
  $f^{-n}$ sends the germ of $A$ at $a$ to $A$ (Claim 3). It transforms the germs of 
  the intervals $I_s\subset A$, $s=1,2$, into germs of two analytic curves intersected 
  transversely at $y_{-n}\neq a$. Thus, $A$ contains the latter germs  and hence, 
  cannot be a finite union of line intervals through $a$, - a contradiction 
  to Claim 1. Proposition \ref{pau} is proved.
  \end{proof}
  
\begin{proof} {\bf of Proposition \ref{afia}.} The line $\hat A$ is 
$\hf$- invariant (Claim 2). Let $(c_{-n},c_{-n-1},\dots)$ be a tail of 
$\hc$, thus, $\hf^{-n}(\hc)\in L(\ha,U)$. By definition, 
$\psi_{\hc}(a)=\hc$. Let $f^{-n}$ be the 
inverse branch germ at $a$ sending $a$ to $c_{-n}\in U$. 
The germ at $a$ of the mapping $\hf^{-n}\circ\psi_{\hc}$ is the composition of two 
mappings: the above germ $f^{-n}$ and $\pi_0^{-1}:U\to L(\ha,U)$. 
The former sends the germ of $A$ at $a$ to $A$ (Claim 3). The latter sends $A$ to 
$\hat A$ by definition. Therefore, $\hf^{-n}\circ\psi_{\hc}(A)\subset\hat A$, 
and hence, $\psi_{\hc}(A)\subset \hat A$. This proves Proposition \ref{afia}. 
\end{proof} 
   
  \section{Accumulation to the horospheres over repellers. Proof of Theorem 
  \ref{taccum}} 
  \subsection{The plan of the proof of Theorem \ref{taccum}} 
Let $f:\oc\to\oc$ be a critically-nonrecurrent rational mapping, $L\subset\caf$ be a leaf 
of the corresponding affine lamination whose projection $p(L)\subset\caf^n$ 
does not lie in a leaf associated to a parabolic 
periodic point. Let $H\subset\chf$ be the corresponding hyperbolic leaf. 
We show that there exists a repelling periodic point $a\in\oc$ of $f$ that is 
not branch-exceptional and 
such that for each horosphere $S\subset H$ the union of its images 
under all the forward and backward iterations of $\hf$ accumulates to some 
point in $H(\ha)$, and hence, to the  horosphere passing through this point. 
This will prove Theorem \ref{taccum}. 

As in the proof of Theorem \ref{trep}, it suffices to prove the above statement 
for just one horosphere $S\subset H$ 
(let us fix it). First we construct appropriate sequence $n_k\in\mathbb N$, 
$n_k\to\infty$, as $k\to\infty$, a  disk $V\subset\oc$ intersecting the Julia set 
(we fix a Hermitian 
Euclidean metric on $V$), appropriate local leaves $L_k\subset\hf^{-n_k}(L)$ 
over $V$ such that 
\begin{equation}\text{for any}\  k\in\nn \ \text{the local leaf} \ 
L_k \  \text{is univalent over} \  V. \label{lkuniv}\end{equation}
For any $w\in V$ we set 
\begin{equation}\hw^k=(\pi_0|_{L_k})^{-1}(w), S^k=\hf^{-n_k}(S), \ 
h_k(w)= \ \text{the height of} \ S^k \ \text{over} \ \hw^k.\label{skwk}\end{equation}
We show that 
\begin{equation}h_k\to-\infty \ \text{uniformly on compact sets in} \ V, \ 
\text{as} \ k\to+\infty.\label{hkinf}\end{equation} 

\begin{proposition} \label{akhk} Let $V$, $L_k$, $S^k$ be as above and satisfy 
(\ref{lkuniv}) and (\ref{hkinf}). Let 
$a\in V$ be an arbitrary  repelling periodic point, $\ha$ be its periodic orbit, 
 $\ha^k$ be as in (\ref{skwk}), 
$\alpha^k=(\ha^k,\theta_k)\in S^k$ be the point of the horosphere $S^k$ 
over $\ha^k$. 
There exist a $h\in\rr$ and  a sequence $l_k\in\nn$, $l_k\to\infty$, as $k\to\infty$, 
such that, after passing to a subsequence, one has 
\begin{equation}\hf^{l_k}(\a^k) \to(\ha,h)\in H(\ha), \ \text{as} \ 
k\to\infty.\label{akconv}
\end{equation} 
\end{proposition}
\begin{proof} 
Let $s\in\nn$ be the period of $a$, $\lambda=|(f^s)'(a)|>1$. 
Then for any $m\in\mathbb N$ one has 
$$\hf^{ms}(\a^k)=(\hf^{ms}(\ha^k),h_{m,k}), \ h_{m,k}=\theta_k+m\ln\lambda. \ \text{Set} \ m_k=
-[\frac{\theta_k}{\ln\lambda}].$$
By construction, the sequence $h_{m_k,k}$ is bounded: 
$0\leq h_{m_k,k}\leq\ln\lambda$.  
Passing to a subsequence one can achieve that it converges to some $h\in\rr$. 
Set $l_k=m_ks$. One has $\theta_k\to-\infty$, by  
(\ref{hkinf}); thus, $m_k,l_k\to+\infty$. Hence, $\hf^{l_k}(\ha^k)\to\ha$, by (\ref{lkuniv}) 
and Corollary \ref{cconva}. 
Then $\hf^{l_k}(\a^k)=(\hf^{l_k}(\ha^k),h_{m_k,k})\to\a=(\ha,h)\in H(\ha)$, by 
Proposition \ref{convh}. This proves (\ref{akconv}). 
\end{proof}

\begin{proof} {\bf of Theorem \ref{taccum} modulo (\ref{lkuniv}) and (\ref{hkinf}).} 
Fix an arbitrary repelling periodic point $a\in V$ (these points are dense in the Julia set). 
 It is not branch-exceptional: there are infinitely many univalent local leaves $L_k$ 
 over $V$, by (\ref{lkuniv}) and (\ref{hkinf}).  
 Let $l_k$, $\alpha^k$, $h$ be as in Proposition \ref{akhk}. 
Then the horospheres $\hf^{l_k}(S^k)$ accumulate to the 
horosphere through $(\ha,h)$ (Proposition \ref{akhk}). This proves 
Theorem \ref{taccum}.
\end{proof}

In the proof of (\ref{lkuniv}) and  (\ref{hkinf}) we use the following 
well-known theorem.  

\begin{theorem} \label{mane} (Ma\~ne, \cite{manie}, theorem II). 
Let $f$ be  a rational function, 
$z\in\overline{\cc}$ be a point that is neither an attractive, nor a parabolic periodic point 
and does not belong to the $\omega$- limit set of a 
recurrent critical point. Then for any $\var>0$ there exists a neighborhood $V=V(z)$ 
such that for any $m\in\nn$ each connected component of the preimage $f^{-m}(V)$ 
has diameter less than $\var$.
\end{theorem}

First we construct $V$, $n_k$ and $L_k$  in the case, 
when $L\subset\caf^l$ (the next subsection). The construction in 
the opposite case, which is given in 3.3, is similar but  slightly 
more technical.

\subsection{Case when $L\subset\caf^l$} 
Fix a $\hx\in L$,  $\pi_0(\hx)\in J=J(f)$. 
The above-mentioned  construction is based on the following 

\begin{lemma} \label{univleaf} Let $f$ be a critically-nonrecurrent rational function, 
$L\subset\caf^l$ be a leaf that is not associated to a parabolic periodic point. 
Let $\hx$ be as above. There exist a sequence $n_k\to+\infty$, a point 
$b\in J(f)$ (that is not a parabolic periodic point) and  a neighborhood 
$V=V(b)\subset\oc$ such that $ x_{-n_k}\to b$ and for any $k\in\mathbb N$ the local 
leaf $L_k=L(\hf^{-n_k}(\hx),V)$ is  univalent over $V$.
\end{lemma}

\begin{proof} 
Recall that $f$ is critically-nonrecurrent. This implies that there is a natural ordering 
of some pairs of critical points of $f$ that lie in the Julia set of $f$. Namely, given two critical points 
$c_1\neq c_2\in J$, we say that 
$c_1>c_2$, if either $c_2$ is the image of $c_1$ under some positive iteration of $f$, or $c_2$ belongs 
to the $\omega$- limit set of $c_1$. 
This ordering has the transitivity property: if $c_1>c_2$ and $c_2>c_3$, then $c_1>c_3$ (critical 
nonrecurrence and absence of periodic critical points in the Julia set). 

Let us consider two different cases.

Case 1): the limit set of $\{ x_{-n}\}|_{n\in\mathbb N}$ contains no critical point of $f$. 
Choose 
a subsequence $n_1<n_2<\dots$ so that $ x_{-n_k}$ converge to some point $b\in\oc$ that  is not a parabolic periodic point. It is possible, since the Riemann sphere is compact and the leaf $L=L(\hx)$ is not a 
leaf associated to a parabolic periodic point. Then for any $\var>0$ there exists 
a neighborhood $V=V(b)$ such that for any $m\in\nn\cup0$ each connected 
component of the preimage $f^{-m}(V)$ has diameter less than $\var$ (by Ma\~ne's 
Theorem). Fix these $\var$ and $V$ so that  
\begin{equation}\var<dist( x_{-n},c) \ \text{for any critical point} \ 
c \ \text{of} \ f \ \text{and any} \ n\in\nn \ \text{large enough}.\label{var}\end{equation}
One can achieve this by assumption: $x_{-n}$ accumulate to no critical point. 
Without loss of generality we assume that $x_{-n_k}\in V$ for any $k$ and 
inequality (\ref{var}) holds true for all $n\geq n_1$. 
One can achieve this by 
removing all ``too small'' $n_k$ by the convergence $x_{-n_k}\to b\in V$. For any 
$m\in\mathbb N$ set 
\begin{equation}
V_{-m,k}=\text{ the connected component containing } x_{-m-n_k} \text{ of } 
f^{-m}(V).\label{vmk}\end{equation}
Each $V_{-m,k}$ contains no critical point by (\ref{var}) and 
since $diam V_{-m,k}<\var$. This implies that  
each local leaf $L_k=L(\hf^{-n_k}(\hx),V)$ is  univalent over $V$. 

Case 2): the sequence $ x_{-n}$ accumulates to some critical point of $f$. 
The limit critical points lie in the Julia set, as does $ x_0$. Let $b$ be a maximal 
 limit critical point. Fix a $\var>0$ satisfying (\ref{var}) for any critical point $c>b$ and 
 a neighborhood $V=V(b)$ satisfying the statement of Ma\~ne's Theorem. Without loss 
 of generality we consider that if $V$ contains an image $f^n(c')$, $n\in\mathbb N$, 
 of a critical point $c'$, then $c'>b$ (hence, $c'$ is not a limit critical point). One can 
 achieve this by shrinking $V$, by the definition of the order on critical points. Fix an 
 arbitrary  subsequence $n_k$ so that 
$ x_{-n_k}\to b$, as $k\to\infty$, and $x_{-n_k}\in V$ for all $k$. Then each 
above local leaf 
 $L_k$ is  univalent over $V$. Indeed,  for any $m\in\nn$ the domain 
 $V_{-m,k}$ from (\ref{vmk}) contains no critical point 
 $c>b$ (by the choice of $\var$ and Ma\~ne's Theorem). It contains no other critical 
 point by the choice of $V$. Lemma \ref{univleaf} is proved.
\end{proof}

Let $n_k$, $b$, $V$ and $L_k$ be as in the above lemma. Statement 
(\ref{lkuniv}) follows from the lemma. In the proof of (\ref{hkinf}) we use the 
following

\begin{proposition} \label{hnorm} For any rational function $f$, any simply connected 
domain $V\subset\oc$ (we fix a Hermitian Euclidean metric on it), any 
 family of univalent local leaves $L_k\subset\caf$ over $V$, any family of horospheres 
 $S^k$ in the corresponding hyperbolic leaves in $\chf$, the  
 height functions $h_k(w)$, see (\ref{skwk}), are equicontinuous on compact 
 sets.
 \end{proposition}
 
 \begin{proof} Without loss of generality we consider that $V$ is 
 unit disk and $b=0$. Each local leaf $L_k$ is affine isomorphic to a domain 
 $V^k\subset\cc$. Indeed, the ambient leaf is isomorphic to the quotient of $\cc$ 
 by a discrete group of Euclidean isometries. The local leaf $L_k$ 
 is simply connected, as is $V$, and contains no singularities  
 (Definition \ref{univdefn}). By definition, one has 
 $$h_k=-\ln|\tau_k'|+const(k), \text{ where }  \tau_k=(\pi_0|_{L_k})^{-1}:
 V\to L_k=V^k\subset\cc.$$ 
 The functions $\tau_k$ are univalent and can be normalized so that 
 $\tau_k(0)=0$, $\tau'_k(0)=1$, by applying affine transformations in the image. 
The latter logarithms remain unchanged (up to additive constants).  
Thus normalized univalent functions form a normal family. This implies the 
 equicontinuity of the heights $h_k$ on compact sets and proves the proposition. 
\end{proof}

\begin{proof} {\bf of (\ref{hkinf}).} Set $l_k=n_k-n_1$. By definition, 
$x_{-n_1}=f^{l_k}(x_{-n_k})$. 
The heights of $S^k$ over the points $\hf^{-n_k}(\hx)$ tend to minus infinity. 
Indeed, they are equal to $\ln|(f^{-l_k})'( x_{-n_1})|$ 
plus the height of $S^1$ over $\hf^{-n_1}(\hx)$ (with respect to the 
standard Euclidean metric in a chart near $ x_{-n_1}$). The latter logarithm tends to 
$-\infty$, since $(f^{-l_k})'( x_{-n_1})\to0$ (the Shrinking Lemma). Thus, 
$h_k(x_{-n_k})\to-\infty$, $x_{-n_k}\to b=0$, 
by construction. This together with Proposition \ref{hnorm} proves 
 (\ref{hkinf}) and finishes the proof of Theorem 
\ref{taccum} in the case, when $L\subset\caf^l$.
\end{proof}

\subsection{Case when $L\subset\caf\setminus\caf^l$} 
 In this case the function $p:L\to p(L)\subset\caf^n$ is not an affine isomorphism. 
 Lemma \ref{univleaf}  and hence, the above construction do not apply 
 literally, and we have to modify our arguments. We fix  a 
$\hx\in L, \ \pi_0(\hx)\in J$,  and a sequence  $n_k\in\nn$, set 
$$\hy=p(\hx)\in\caf^n=\caf^l, 
\text{ such that } n_k\to\infty, \ y_{-n_k}\to b\in\oc, \text{ as } k\to\infty,$$
$b$ is not a parabolic periodic point. The existence of these $n_k$ and $b$ follows 
from the condition of Theorem \ref{taccum}, which says that the projection 
$p(L)\subset\caf^n$ is not contained in a leaf 
associated to a parabolic periodic point.  For any critical 
point $c$ of $f$ denote $d_c$ the local degree of $f$ at $c$: the multiplicity of $c$ as a 
preimage of its critical value. Set 
\begin{equation}
d=\prod_{c\in J} d_c: \text{ the product being taken through all the critical points of } f 
\text{ in } J.
\label{defd}\end{equation}
We show that passing to a subsequence, one can achieve that 
there exists a disk $U\subset\oc$ centered at $b$ such that each uniformized 
local leaf $L(\hf^{-n_k}(\hx),U)$ 
(see the next definition) is a branched cover over $U$ 
of one and the same degree $\nu\leq d$, and each local leaf $L(\hf^{-n_k}(\hy),U)$ 
is univalent (Lemmas \ref{univleaf}, \ref{lboundeg} and Corollary \ref{cornk}). 
We prove that the ramification points of the above uniformized local leaves  
tend to $b$ (Proposition \ref{gconv}, which  provides a stronger statement). 
We parametrize each  uniformized local leaf by unit disk $D_1$ equipped 
with the standard Euclidean metric. We show that the corresponding heights of the 
uniformized horospheres $S^k$ tend to $-\infty$ uniformly on compact sets in $D_1$ 
(Proposition \ref{wthkinf}). Afterwards we fix a disk $V$ such that $\overline V\subset 
U\setminus b$, and for any $k$ large enough we fix a local leaf 
$L_k\subset L(\hf^{-n_k}(\hx),U)$ over $V$. The convergence of the ramification points 
implies the univalence of $L_k$ for large $k$. The  height convergence 
statement (\ref{hkinf}) for thus constructed $L_k$ will be then deduced from 
Propositions \ref{gconv}  and \ref{wthkinf}. 

\begin{definition} \label{unileaf}
Let $L\subset\caf$ be a leaf, $B\subset L$ be an open domain. 
Recall that $L$ is isomorphic to either $\cc$, or its quotient by a discrete group of 
Euclidean isometries. 
 Let $\Pi:\cc\to L$  be the quotient projection, $\wt B\subset\cc$ be a connected 
component of $\Pi^{-1}(B)$. 
The domain $\wt B$ and the projection $\Pi:\wt B\to B$ are called the {\it affine 
uniformization} of $B$. If $B$ is a (local) leaf, then $\wt B$ is called the 
{\it uniformized (local) leaf}, and the composition $\pi_0\circ\Pi:\wt B\to\pi_0(B)$ is called 
the {\it uniformizing (local) leaf projection}. Let $H$ be the hyperbolic leaf in $\chf$ 
corresponding to $L$. The projection $\Pi$ extends up to a quotient projection 
$\Pi:\htr\to H$ called the  uniformization of $H$. 
For any horosphere $S\subset H$ its preimage $\wt S=\Pi^{-1}(S)\subset\htr$ is a 
horosphere called the {\it uniformizing horosphere} of $S$. 
\end{definition}

\begin{remark} Let $L\subset\caf$ be a leaf, $\hx\in L$, 
 \begin{equation}
 \Phi_{\hx}(t)=(\phi_{0,\hx},\phi_{-1,\hx},\dots)(t), \ \Phi_{\hx}(0)=p(\hx)\in\caf^n: \ 
 \phi_{-j,\hx}(0)=\pi_{-j}(\hx),\label{merosec}\end{equation}
 be the meromorphic function sequence (\ref{funcsec}) defining $\hx\in\caf$. For 
 each $w\in\cc$ the meromorphic function sequence $\Phi_{\hx}(w+t)$ represents 
 a point of $L$ denoted $\hw$.
 The mapping $\cc\to L$: $w\mapsto\hw$ is an affine uniformization,  
 and $\phi_{0,\hx}:\cc\to\overline{\cc}$ is the uniformizing projection.
A local leaf  is univalent, if and only if so is the uniformizing local leaf 
projection. 
\end{remark}

\begin{lemma} \label{lboundeg} Let $f$ be a critically-nonrecurrent rational 
function,  $d$ be as in (\ref{defd}). For any $b\in J=J(f)$ that is not a parabolic 
periodic point there exists a disk $U\subset\cc$ centered at $b$ 
such that for any $\hb\in\caf$ with $\pi_0(\hb)=b$ the corresponding uniformized 
local leaf over $U$ is simply connected and is a branched cover over $U$ of 
degree no greater than $d$.
\end{lemma}

\begin{proof} By critical nonrecurrence, there exists a $\var>0$ (let us fix it) such that 
\begin{equation} dist(f^n(c),c)>\var \ \text{for any critical point} \ c\in J=J(f) \ 
 \text{of} \ f, \ n\in\nn,\label{distcr1}\end{equation}
\begin{equation} \ dist(c,c')>\var
  \ \text{for any two distinct (may be multiple) critical points} \ 
  c, c'\in J.\label{distcr2}\end{equation}
Fix an arbitrary disk $U\subset\oc$ centered at $b$ and disjoint 
from the forward orbits of the critical points in the Fatou set,  
such that for any $n\in\nn\cup 0$ 
each connected component of the preimage $f^{-n}(U)$ has diameter less than 
$\var$ (Ma\~ne's Theorem). Let us prove the lemma for this $U$. To do this, we 
use the following 

\medskip

{\bf Claim 1.} {\it Let $V,W\subset\oc$ be connected domains, $W$ be 
homeomorphic to disk. Let  
$f:V\to W$ be a finite degree branched covering with a unique critical 
point $c$. Let $d$ be the local degree of $c$. 
Then $V$ is simply connected and $deg f=d$.}

\medskip

Claim 1 follows immediately from the Riemann-Hurwitz Formula.

Case 1): $\hb\in\caf^l$.  For any $n$ let $U_{-n}$ be the above connected component 
that contains $b_{-n}$. Each $U_{-n}$ contains at most one (may be multiple) 
critical point of $f$, 
by the inequalities $diam(U_{-n})<\var$ and (\ref{distcr2}). A critical point $c$ 
cannot be contained in two domains $U_{-n}$ and $U_{-n-m}$: otherwise 
$c,f^m(c)\in U_{-n}$, which is impossible by the diameter bound and (\ref{distcr1}). 
Thus, each mapping $f:U_{-n}\to U_{-n+1}$ is either a
conformal diffeomorphism, or a branched covering with a unique critical point 
$c_n\in J$. The domains $U_n$ are simply connected, and the degrees of the 
latter coverings are equal respectively to 1 and $d_{c_n}$ (Claim 1 and induction 
in $n$).  The number of the values of $n$ corresponding to branchings is no greater 
than the number of the critical points in the Julia set. 
Let $m$ be the maximal one of these $n$. Then $\pi_{-m}:L(\hb,U)\to U_{-m}$ is 
a diffeomorphism (recall that $\hb\in\caf^l$ by assumption). 
Hence, the local leaf $L(\hb,U)$ is simply connected, as is 
$U_{-m}$,  and $\pi_0=f^m\circ\pi_{-m}:L(\hb,U)\to U$ is a branched cover of 
degree $\prod d_{c_n}\leq d$.  The local leaf $L(\hb,U)$ is isomorphic to its affine uniformization, since 
it is contained in $\caf^l$. This proves the statement of the lemma for  $\hb\in\caf^l$. 
 
 Case 2): $\hb\in\caf\setminus\caf^l$. The set $\caf^l$ is dense in $\caf$ 
 (Proposition \ref{minim}). Fix a  sequence $\hb^n\in\caf^l$, $b^n_0=b$, 
 $\hb^n\to\hb$  in $\caf$. Let $\phi_{0,\hb^n}$, $\phi_{0,\hb}$ be the corresponding 
 functions from (\ref{merosec}) normalized so that  $\phi_{0,\hb^n}\to\phi_{0,\hb}$ 
  uniformly on compact sets in $\cc$. 
 Recall that their values at 0 equal $b$. Set 
 $B^n=L(\hb^n,U)$. Let $\wt B^n$ ($\wt B$) denote the connected component 
 containing $0$ of the preimage $\phi_{0,\hb^n}^{-1}(U)$ (respectively, 
 $\phi_{0,\hb}^{-1}(U)$). Then each mapping 
 $\phi_{0,\hb^n}:\wt B^n\to U$ is the uniformizing 
 projection corresponding to the local leaf $B^n$. Hence, it is a branched cover of 
 degree at most $d$, and $\wt B^n$ is simply connected, as 
was proved above. 

\medskip

{\bf Claim 2.} {\it The domain $\wt B$ is simply connected,  and the  
mapping $\phi_{0,\hb}:\wt B\to U$ is a branched cover of degree at most $d$.}

\begin{proof} The simply connected domains $\wt B^n$ are conformally equivalent 
to the unit disk $D_1$, since $\phi_{0,\hb^n}:\wt B^n\to U$ are nonconstant 
bounded holomorphic functions. Fix some conformal isomorphisms 
$\eta_n:D_1\to\wt B^n$, $\eta_n(0)=0$. Passing to a subsequence, one can achieve 
that $\eta_n$ converge to a nonconstant univalent function $\eta:D_1\to\cc$. 
This follows from the compactness of the space of normalized univalent functions, 
K\"obe $\frac 14$ theorem and 
the uniform boundedness of the distances $dist(\partial\wt B^n,0)$ from above and 
from below. Indeed,  $dist(\partial\wt B^n,0)\to dist(\partial\wt B,0)$, as $n\to\infty$, 
by the convergence $\phi_{0,\hb^n}\to\phi_{0,\hb}$. The mappings 
$g_n=\phi_{0,\hb^n}\circ\eta_n:D_1\to U$ are Bl\"aschke products of degree at most $d$, 
since they are branched covers of degree at most $d$, as are 
$\phi_{0,\hb^n}:\wt B^n\to U$,  and $D_1$, $U$ are disks. One has  
$g_n\to g=\phi_{0,\hb}\circ\eta\not\equiv const$ uniformly on compact sets. 
The limit $g$ is a Bl\"aschke product of degree at most $d$, as are $g_n$.  Therefore, 
the mapping $\phi_{0,\hb}:\eta(D_1)\to U$ is a branched cover of degree at most 
$d$. The domain $\eta(D_1)$ is simply connected (univalence) and coincides 
with $\wt B$ by construction. The claim is proved. 
\end{proof}

 The mapping from Claim 2 is the uniformizing projection corresponding to the 
 local leaf $L(\hb,U)$. This proves the lemma.
 \end{proof}
 
\begin{corollary} \label{cornk} Let $\hx\in\caf$,  $\hy=p(\hx)\in\caf^l$, 
$y_0=\pi_0(\hx)\in J$,  and $L(\hy)$ be not a 
leaf associated to a parabolic periodic point. Let $b\in J$ be a limit point of 
the sequence $y_0,y_{-1},\dots$ that is not a parabolic periodic point. There exist 
a disk $U\subset\oc$ centered at $b$, a $\nu\in\nn$, $\nu\leq d$, 
and a sequence $n_k\in\nn$ satisfying the following statements: 
\begin{equation}y_{-n_k}=\pi_{-n_k}(\hx)\to b, \text{ as } k\to\infty;\label{convnk}
\end{equation}
\begin{equation} \text{ the local leaves } L(\hf^{-n_k}(\hy),U) \text{ are univalent over } 
U;\label{univy}\end{equation}
the uniformizing projections corresponding to the local leaves $L(\hf^{-n_k}(\hx),U)$ 
are branched covers over $U$ of degree $\nu$. 
\end{corollary}
 
 \begin{proof} There exist a sequence $n_k$ and a disk $U_1$ centered at $b$ 
 that satisfy (\ref{convnk}) and (\ref{univy}) (Lemma \ref{univleaf}). There exists another 
  disk $U_2$ centered at $b$ and satisfying the statements of Lemma \ref{lboundeg}. Let $U$ be the 
  smallest one of the disks $U_1$ and $U_2$. The local leaves 
  $L(\hf^{-n_k}(\hy),U)$ remain univalent after passing to a smaller disk.  
  The uniformized local leaves $L(\hf^{-n_k}(\hx),U)$ remain branched covers 
  over $U$ of degrees $\nu_k\leq d$. 
 Passing to a subsequence one can achieve that the degrees $\nu_k$ are the same. 
 This proves the corollary.   
  \end{proof} 
  
 We fix a $\hx\in L$, set $\hy=p(\hx)$, 
 and $b$, $n_k$, $\nu$, $U$, as in Corollary \ref{cornk}. Without loss of 
 generality we consider that there is a bigger disk $U'\Supset U$ 
 over which the local leaves $L(\hf^{-n_k}(\hy),U')$ are univalent 
 (shrinking $U$ and passing to a subsequence, as in the proof of the corollary). 
 For any $k\geq2$ set  
 $$l_k=n_k-n_1, \ f^{-l_k}:U'\to\oc \text{ the inverse branches such that } 
 f^{-l_k}(y_{-n_1})=y_{-n_k}.$$
 These branches are holomorphic, by the  univalence of the local leaves 
 $L(\hf^{-n_k}(\hy),U')$. By Shrinking Lemma and (\ref{convnk}), one has  
 \begin{equation}f^{-l_k}|_U\to b \text{ uniformly on } U. \text{ We consider that } 
 f^{-l_k}(U)\subset U\label{uincl}\end{equation}
 for all $k$, passing to an appropriate  subsequence. 
 Without loss of generality everywhere below we consider that 
 $$U=D_1, \ b=0. \text{ Set } \hb^k=\hf^{-n_k}(\hx).$$
 For any $k$ let $\wt B^k\subset\cc$ denote the connected component containing 0 of the 
 preimage $\phi_{0,\hb^k}^{-1}(U)$, i.e., the uniformized local leaf 
 $L(\hb^k,U)$. The domains $\wt B^k$ are simply connected,  
 and $\phi_{0,\hb^k}:\wt B^k\to U$ are branched covers of degree $\nu$ 
 (Corollary \ref{cornk}). Fix conformal isomorphisms $\eta_k:D_1\to\wt B^k$, 
 $\eta_k(0)=0$. Set 
 \begin{equation}
 g_k=\phi_{0,\hb^k}\circ\eta_k:D_1\to U=D_1, \label{gk}\end{equation}
 which are branched covers of degree $\nu$, and hence, Bl\"aschke products,  
 \begin{equation} g_k(0)=0, \ g_k(1)=1,\label{ngn}\end{equation} 
 after appropriate normalization of $\eta_k$ by rotation in the source.   
 \begin{proposition} \label{gconv} One has $g_k(z)\to g(z)=z^{\nu}$ uniformly on 
 $\overline{D_1}$.
 \end{proposition}
 
 \begin{proof} Set  
 $$C_k=\{ \text{the critical values of } \phi_{0,\hb^k} \} =
 \{ \text{the critical values of } g_k\} \subset U=D_1.$$
 Each $C_k$ is a set of $\nu-1$ points, some of them may coincide. One has 
 \begin{equation} C_k\to b=0, \text{ as } k\to\infty,\label{ck}\end{equation}
by the convergence $f^{-l_k}|_U\to 0$ and since $f^{-l_k}(C_1)=C_k$. Indeed, 
$f^{-l_k}(C_1)\subset C_k$, since $\phi_{0,\hb^k}=f^{-l_k}\circ\phi_{0,\hb^1}$ up to 
$\cc^*$- action in the source, $f^{-l_k}$ is holomorphic on $U$ and 
$f^{-l_k}(U)\subset U$, see (\ref{uincl}). Therefore, $f^{-l_k}(C_1)=C_k$, since 
the total multiplicities of the critical values in $C_1$ and $C_k$ are equal. 
Passing to a subsequence one can achieve that $g_k$ converge uniformly on compact 
sets and $C_k\subset D_{\frac12}$. 
Let $g$ denote the limit of $g_k$. Let us show that $g(z)=z^{\nu}$. To do this, 
we consider the auxiliary annulus 
$$A=D_1\setminus\overline{D_{\frac12}}, 
\text{ set } \wt A_k=g_k^{-1}(A), \ \wt A=g^{-1}(A).$$
One has $C_k\cap A=\emptyset$ for all $k$. Hence, each mapping 
$g_k:\wt A_k\to A$ is a nonramified covering of degree $\nu$, and $\wt A_k$ is an 
annulus adjacent to $\partial D_1$, $mod(\wt A_k)=\frac1{\nu}mod(A)$. Set 
\begin{equation}W_k=D_1\setminus\wt A_k=g_k^{-1}(\overline{D_{\frac12}}).
\label{adincl}\end{equation}
There exists a $0<r<1$ such that 
\begin{equation} W_k\Subset D_r \text{ for all } k.\label{star2}\end{equation}
 Indeed, $0\in g_k^{-1}(0)\subset W_k$, see (\ref{adincl}). Therefore,  the  
 contrary to (\ref{star2}) would imply that the complement $W_k$ to 
 $\wt A_k$ contains both the origin  and points arbitrarily close to $\partial D_1$. 
 Hence, the moduli $mod(\wt A_k)$ are arbitrarily small, - a contradiction to the equality 
 $mod(\wt A_k)=\frac1{\nu}mod(A)$. 

The limit $g$ is not a constant. Indeed,  otherwise, $g_k\to0=g_k(0)$ and 
$g_k(D_r)\subset D_{\frac12}$ for large $k$, hence, $W_k\supset D_r$, see 
(\ref{adincl}), - a contradiction to (\ref{star2}). 
Therefore, $g$ is a Bl\"aschke product, as are $g_k$. 
The critical points of $g_k$ lie in  $W_k\subset D_r$, since $C_k\subset D_{\frac12}$ 
and by (\ref{adincl}). 
Therefore, passing to a subsequence one can achieve that they converge, 
and their limits are critical points of $g$. Thus, $g$ is a branched cover of degree $\nu$  
with a unique critical value $0$, see (\ref{ck}). In particular, the convergence $g_k\to g$ is 
uniform on $\overline{D_1}$, since the latter is a convergence of Bl\"aschke products  
of a given degree to a Bl\"aschke product of the same degree. 
The preimage $g^{-1}(D_1\setminus 0)$ is 
a finite regular cover over a punctured disk, and hence, conformally equivalent to 
a punctured disk. It is adjacent to $\partial D_1$ and does not contain 0. Hence, 
it is $D_1\setminus0$, 
$g^{-1}(0)=0$ and $g(z)=cz^{\nu}$. Now $c=1$, by (\ref{ngn}). This proves the 
proposition. 
\end{proof}

Fix a horosphere $S\subset H$. Set 
$$S^k=\hf^{-n_k}(S), \ \wt S_k = \text{ the uniformizing horosphere of } S^k,$$
see Definition \ref{unileaf}. Consider the pushforward under 
$\eta_k:D_1\to\wt B^k$ of the standard Euclidean metric on $D_1$. This is an 
Euclidean metric on $\wt B^k$. For any $w\in D_1$ set
$$\wt h_k(w)=\text{ the height of } \wt S_k \text{ over the point } \eta_k(w) 
\text{ in the above metric.}$$

\begin{proposition} \label{wthkinf} One has $\wt h_k\to-\infty$ uniformly on compact 
sets in $D_1$, as $k\to\infty$. 
\end{proposition}

\begin{proof} The mapping $\hf^{-l_k}$ sends the local leaf $L(\hf^{-n_1}(\hx),U)$ to 
$L(\hf^{-n_k}(\hx),U)$, by definition and (\ref{uincl}). Its lifting 
to the affine uniformization is 
an affine mapping $\wt B^1\to\wt B^k$, which will be also denoted $\hf^{-l_k}$. Set 
$$\chi_k=\eta_k^{-1}\circ\hf^{-l_k}\circ\eta_1:D_1\to D_1. 
\text{ One has }$$
\begin{equation} \wt h_k(\chi_k(w))=\wt h_1(w)+\log|\chi'_k(w)|,\label{hkwt}\end{equation}
\begin{equation}\chi_k\to0 \text{ uniformly on compact sets in } D_1, 
\text{ as } k\to\infty,\label{chilim}
\end{equation}
by definition, (\ref{uincl}), Proposition \ref{gconv} and since 
$\chi_k$ coincides with an analytic branch of the composition 
$g_k^{-1}\circ f^{-l_k}\circ g_1:D_1\to D_1$. Hence, $\chi_k'\to0$. 
This together with (\ref{hkwt}) 
implies that $\wt h_k(\chi_k(0))\to-\infty$. Hence, 
$\wt h_k\to-\infty$ uniformly on compact sets, by (\ref{chilim}) and equicontinuity (Proposition \ref{hnorm}). This proves Proposition \ref{wthkinf}. 
\end{proof}

\begin{proof} {\bf of Theorem \ref{taccum}.} Let $\hx$, $n_k$, $U$, $b$ be the same, 
as in Corollary \ref{cornk}. Fix an arbitrary disk $V$ such that 
$\overline V\subset U\setminus b$. The uniformizing projections corresponding to 
the local leaves $L(\hf^{-n_k}(\hx),U)$ are branched covers over $U$. 
Without loss of generality we consider that 
$\overline V$ contains no their critical values, passing to $k$ large enough, see (\ref{ck}). 
For any $k$ fix a local leaf $L_k\subset L(\hf^{-n_k}(\hx),U)$ over $V$. It is univalent by 
construction. This proves statement (\ref{lkuniv}). Let us prove (\ref{hkinf}). Its 
stronger version says that the heights of the horospheres $\wt S_k$ over $L_k$
 in the standard Euclidean metric  of the disk 
$V$ tend to $-\infty$ uniformly on $V$. The latter heights 
differ from the heights $\wt h_k$ by the logarithm of  the ratio of the standard Euclidean 
metric on the disk $U$ and the pushforward under $g_k$ of the Euclidean metric on 
$D_1$. The latter ratio  is uniformly bounded from above and from below on 
$\overline V$, and there 
exists a $0<r<1$ such that $g_k^{-1}(\overline V)\subset \overline D_r$ for all $k$ 
(Proposition \ref{gconv}). This together with Proposition \ref{wthkinf} proves (\ref{hkinf}).  
Thus, statements (\ref{lkuniv}) and (\ref{hkinf}) are proved for the above $L_k$ and $V$. 
This together with the discussion in 3.1 proves Theorem \ref{taccum}.
\end{proof}
  
 \section{Closeness of the horospheres associated to the parabolic periodic points}
Here we prove Theorem \ref{taccum2}.
Let $f$ be a rational function with a parabolic periodic point 
$a\in\oc$. Without loss of generality we consider that 
$a$ is fixed, passing to appropriate iteration of $f$. Let $L_a\subset\caf$ 
be the leaf associated to $a$ of the affine lamination, $H_a\subset\chf$ be the corresponding hyperbolic 
leaf.  In the proof of Theorem \ref{taccum2} we use the following proposition. 

\begin{proposition} \label{horoinv} 
Let $f$, $a$, $H_a$ be as above. Then each horosphere in $H_a$ is invariant 
under the mapping $\hf$. 
\end{proposition}
\begin{proof} The Fatou coordinate is affine on the leaf $L_a$, and $\hf$ acts by unit translation there. 
Hence, it preserves a Euclidean metric on $L_a$. Therefore, the lifting of $\hf$ to 
$H_a$ preserves the corresponding height. This implies the proposition. 
\end{proof}

Fix a horosphere $S\subset H_a$. We show that $S$ is closed in $\chf$ 
and does not accumulate 
to itself. This together with its  invariance (Proposition \ref{horoinv}) implies 
Theorem \ref{taccum2}. 

\def\hu{\hat u}

We prove the previous statement by contradiction. Suppose the contrary: then $S$ 
accumulates to a horosphere $S'\subset\chf$. Let $H\subset\chf$ be the leaf 
containing $S'$, $L\subset\caf$ be the corresponding affine leaf. Fix a 
disk $U\subset\oc$, a $(\hb,h)\in S'$, a sequence $(\hb^k,h_k)
\in S$, $\hb\in L$, $\hb^k\in L_a$, such that 
\begin{equation}
\pi_0(\hb) \in U, \pi_0(\hb)\neq a,\ \text{and the local leaf} \ 
L(\hb,U)\  \text{is univalent 
over} \ U,\label{unibk}\end{equation}
\begin{equation}
(\hb^k,h_k)\to(\hb,h) \text{ as } k\to\infty, \ \pi_0(\hb^k)=\pi_0(\hb), \ 
\hb^k\neq\hb 
\text{ for all } k.\label{bkhk}\end{equation}
Here the heights of the horospheres are measured with respect to the Euclidean 
metric of $U$. Without loss of generality we consider that $U=D_1$ and 
$\pi_0(\hb)=0$. Fix an arbitrary disk $V$ centered at 0 such that $\overline V
\subset U$. Then for any $k$ large enough the local leaf 
$$\Lambda_k=L(\hb^k,V) \ \text{is univalent over} \ V,$$
by (\ref{unibk}) and the definition of topology on $\caf$. 
Without loss of generality we consider  that all $\hb^k$, and hence, 
$\Lambda_k$, are 
distinct. One can achieve this by passing to a subsequence, by (\ref{bkhk}). 
We show that 
\begin{equation}h_k\to +\infty, \text{ as } k\to\infty,\label{hks}\end{equation}   
  - a contradiction to (\ref{bkhk}). This proves Theorem \ref{taccum2}.  

\def\hv{\hat v}
The inverse germ $f^{-1}$ fixing $a$ will be called the {\it parabolic germ}. 
For the proof of (\ref{hks}) we fix a $r>0$, set $\Delta=D_{2r}(a)$, such that 
\begin{equation}\text{ the parabolic inverse germ } f^{-1} \text{ is holomorphic 
on } \Delta\cup f(\Delta) \text{ and } f(\overline{D_r(a)})\subset\Delta; \label{rd1}\end{equation}
\begin{equation} \text{if } z\in\overline{D_r(a)}, 
 \text{ then the parabolic branch } f^{-1}|_{\Delta} \text{ sends } f(z)\in\Delta 
 \text{ to } z;
 \label{rd2}\end{equation}
 \begin{equation} \text{each infinite backward orbit of } f 
 \text{ contained in } \overline{D_r(a)} 
\text{ converges to } a.\label{rd3}\end{equation}
 Statements (\ref{rd1})-(\ref{rd3}) hold true, whenever $r$ is small enough. 
 For any $k$ let  $n_k\in\nn$ 
be the minimal number such that $b^k_{-j}\in\overline{D_r(a)}$, whenever  
$j\geq n_k$. The number $n_k$ exists, since $\hb^k\in L_a$ and by Proposition 
\ref{repell}.  Set 
$$\hv^k=\hf^{-n_k}(\hb^k).$$
Passing to a subsequence, one can achieve that $v^k_0=b^k_{-n_k}$ converge 
to some point $x_0\in\overline{D_r(a)}$. Then $v^k_{-j}\to x_{-j}\in
\overline{D_r(a)}$ for all $j$, and both $\hv^k$ and $\hx=(x_0,x_{-1},\dots)$ are infinite 
backward orbits of $f$ in $\overline{D_r(a)}$, by (\ref{rd2}). 
They converge to $a$ by (\ref{rd3}). They are distinct from the fixed orbit  
of $a$. For the former, $\hv^k$, this follows from definition and the inequality 
$b^k_0=\pi_0(\hb)\neq a$, see (\ref{unibk}). For the latter, $\hx$, 
this holds true, 
since $x_0=\lim_{k\to\infty} b^k_{-n_k}$ and $b^k_{-n_k}$ are bounded away from 
$a$:  either $n_k=0$ and $b^k_{-n_k}=\pi_0(\hb)\neq a$, or 
$b^k_{-n_k+1}=f(b^k_{-n_k})\notin D_r(a)$, 
by definition. The  backward iterations 
$f^{-j}:x_0\mapsto x_{-j}$ are powers of the parabolic inverse branch, by 
(\ref{rd2}). They 
are holomorphic in a neighborhood $W=W(x_0)$ and converge to $a$ uniformly on 
compact subsets in $W$, since the attractive basin of a parabolic fixed point 
is open. Fix this $W$. 
Then the local leaf $L(\hx,W)$ is univalent over $W$, and $\hv^k\to\hx$ along 
this local leaf, by definition and (\ref{rd2}). In particular, 
$L(\hx,W)\subset L_a$, since $\hv^k\in L_a$. 

 The sequence $n_k$ tends to infinity, by definition and since the 
local leaves $\Lambda_k$ are distinct. Fix a metric on $W$. The corresponding height 
of $S=\hf^{-n_k}(S)$ over $\hv^k$ tends to a finite value, namely, to its height over 
$\hx$, as $k\to\infty$ (since $\hv^k\to\hx$ along a fixed local leaf). 
On the other hand, its difference 
with the height $h_k$ of $S$ over $\hb^k$ is equal to $\ln|(f^{-n_k})'(0)|$, which tends to 
$-\infty$ (by the Shrinking Lemma). 
This implies (\ref{hks}) and  proves Theorem \ref{taccum2}.

 \section{Nondense horospheres. Proof of Theorem \ref{tnodense}}
 
Here we give a proof of Theorem \ref{tnodense} (Subsections 5.1 -- 5.4). 
 Its Addendum is proved analogously with small modifications discussed 
  at the end of the paper. 
 
 \subsection{The plan of the proof of Theorem \ref{tnodense}}
  It suffices to prove Theorem \ref{tnodense} for $h=0$, as in the proof of 
 Theorem \ref{trep} in Section 2. 
 
 Let $\var<\frac14$, $\var\neq0,-2$. 
Let $a=a(\var)$ be the fixed point (\ref{aedef}) of $\fe$, $\ha\in\cafe$ be its fixed orbit, 
  $\Pi_a\subset\cafe$ be the corresponding 
  subset from (\ref{univa}). The set $\Pi_a$ is nonempty, since the fixed point 
  $a(\var)$ is not branch-exceptional (Propositions \ref{prex} and 
 \ref{pianep}).  We prove Theorem \ref{tnodense} for 
 \begin{equation}\be=\overline{\{\b(\ha,\hy) \ | \ \hy\in\Pi_a\}}. 
 \text{ Set} \ SB_\var=\cup_{\b\in\be}S_{\ha,\b}\subset \chfe.
 \label{defbe}\end{equation}
 To do this, we show that for every $\var<\frac14$, $\var\neq0,-2$,  
\begin{equation} \be=\{\sum_{m=1}^l\b(\ha,\hy(m)) \ | \ l\in\nn, \ 
\hy(m)\in\Pi_a\};
\label{countsum}\end{equation}
 \begin{equation}  S_{\ha,0} \ \text{accumulates exactly to} \ 
 SB_\var.\label{closhor}
 \end{equation}
 Statement (\ref{closhor}) is proved in 5.4. 
 Statements (\ref{countsum}) is proved below and in 5.2, 5.3. 
 In the proof of (\ref{countsum}) we use the well-known formula for 
 the basic cocycle $\beta(\hx,\hy)$ as an infinite series in $\ln|\fe'(x_{-j})|$ and 
 $\ln|\fe'(y_{-j})|$ (the next proposition and (\ref{bxx})). We show 
 (Lemma \ref{jul} and Corollary \ref{cjul})  that if 
 $\var<0$ ($0<\var<\frac14$), then the module $|\fe'|$ achieves its maximum 
 (minimum) on the Julia set exactly at $\pm a(\var)$. This together with 
 (\ref{bxx}) implies that $\pm\beta(\ha,\hy)>0$ whenever $\pm\var>0$ and 
 $\hy\in\Pi_a$ (Corollary \ref{corsign}). (Moreover, for a given $\var$, the values $\beta(\ha,\hy)$ 
 are bounded away from $0$.) Using this statement, we show (Lemma \ref{ym} and 
 Corollary \ref{corconv}) 
 that for every sequence $\ha^n\in\Pi_a$ with converging values $\beta(\ha,\ha^n)$ 
 the limit of the latter is a sum (\ref{countsum}). This implies (\ref{countsum}). 
 
 \begin{proposition} \label{busfor} Let $f(z)$ be a rational function, $\hx,\hy\in\caf^l$ 
 be a pair of points lying on one and the same leaf, 
 $ x_{-j}, y_{-j}\in\cc$ for all $j$ (the equality $ x_0= y_0$ does not 
 necessarily hold). Let $\b_{\hx}$, $\b_{\hy}$ be the 
 height functions (\ref{defhei}) on the corresponding 
 hyperbolic leaf in $\chf$ defined by the standard Euclidean metric on $\cc$. Then 
 \begin{equation}\b(\hx,\hy)=\b_{\hy}-\b_{\hx}=\sum_{j=1}^{+\infty}(\ln|f'( y_{-j})|-\ln|f'( x_{-j})|).\label{bxx'}
 \end{equation}
 \end{proposition}
 
 The proposition follows from formula (3.23) in \cite{ly}. It implies that 
 for every $\hy\in\Pi_a$ 
 \begin{equation}\b(\ha,\hy)=\sum_{j=1}^{+\infty}(\ln|\fe'(y_{-j})|-\ln|\fe'(a(\var))|)
 \label{bxx}\end{equation} 
   
  \begin{lemma} \label{jul} The Julia set $J=J(\fe)$ is 
  
  - contained in $D_{a(\var)}\cup\pm a(\var)$, $D_{a(\var)}=\{|z|<a(\var)\}$, 
  whenever $\var<0$;
  
  -  contained in the complement $\cc\setminus(\overline{D_{a(\var)}}\setminus\pm a(\var))$, if 
  $0<\var<\frac14$.
  \end{lemma}
  
  Lemma \ref{jul} is proved in the next subsection. 
  \begin{corollary} \label{cjul}  If $\var<0$ ($0<\var<\frac14$), then for every 
  $x\in J(\fe)$ one has 
  $|\fe'(x)|\leq\fe'(a(\var))$ (respectively, $|\fe'(x)|\geq\fe'(a(\var))$). 
  In both cases the equality is achieved  exactly at $x=\pm a(\var)$.
  \end{corollary}
  
  \begin{proof} One has $\fe'(x)=2x$. 
  This together with the lemma implies the corollary.
  \end{proof}
  
  \def\dse{D_{\sigma}(a(\var))}
  \def\ae{a(\var)}
  \def\hae{\ha(\var)}
  \begin{corollary} \label{corsign} 
  If $\var\in(-\infty,0)\setminus\{-2\}$ 
  ($0<\var<\frac14$), then for every $\hy\in\Pi_a$ each 
  corresponding nonzero term of the sum in (\ref{bxx}), and hence, $\b(\ha,\hy)$, 
  is negative (respectively, positive). The zero terms correspond exactly to 
  $y_{-j}=\pm a(\var)$. 
  \end{corollary}
  
  \begin{lemma} \label{ym} Let $\var\in(-\infty,\frac14)\setminus\{-2,0\}$. 
  Let $\ha^n\in\Pi_a$ be a sequence of points such that the values 
  $\beta(\ha,\ha^n)$ converge to a finite limit and $a^n_{-1}=-a(\var)$. 
  Then (passing to a subsequence 
  one can achieve that) there exists a $l\in\nn$ 
  such that for every $n$ there exist indices 
  $$0=\nu_1<\nu_2<\dots<\nu_l, \ \nu_m=\nu_m(\ha^n), \text{ we set } \nu_{l+1}=
  +\infty,$$
  such that each $m=1,\dots,l$ satisfies the following statements:
  \begin{equation}a^n_{-\nu_m(\ha^n)}\to a(\var), \ \text{as} \ n\to\infty, 
  \label{convnu}\end{equation}
  \begin{equation}\nu_{m+1}(\ha^n)-\nu_m(\ha^n)\to+\infty,\label{diffinf}
  \end{equation}
  \begin{equation}\text{for every } j\in\mathbb N 
  \text{ the sequence} \ a^n_{-\nu_m(\ha^n)-j} \ 
  \text{converges}; \ 
  \text{set} \ y_{-j}(m)=\lim_{n\to\infty}a^n_{-\nu_m(\ha^n)-j},
  \label{defam}\end{equation}
  \begin{equation} \hy(m)=(y_0(m),y_{-1}(m),\dots) \ 
  \text{represents a  point from} 
  \ \Pi_a,\label{pia}\end{equation}
 \begin{equation} \Sigma_m=\sum_{j=\nu_m(\ha^n)+1}^{\nu_{m+1}(\ha^n)}
 (\ln|\fe'(a^n_{-j})|-
  \ln|\fe'(a(\var))|)\to\b(\ha,\hy(m)), \ \text{as} \ n\to\infty.\label{sigmam}
  \end{equation} 
  \end{lemma}
  
  \begin{corollary} \label{corconv} In Lemma \ref{ym} one has 
  $\beta(\ha,\ha^n)\to\sum_{m=1}^l\beta(\ha,\hy(m))$, as $n\to\infty$.
  \end{corollary}
  The corollary follows from (\ref{bxx}) and (\ref{sigmam}).
  
  \begin{proof} {\bf of Theorem \ref{tnodense} modulo 
  Lemmas \ref{jul}, \ref{ym} and (\ref{closhor}).} Consider an arbitrary 
  converging sequence of basic cocycle values $\beta(\ha,\ha^n)$, 
  $\ha^n\in\Pi_a$. Without loss of generality we consider that 
  $a^n_{-1}\neq a(\var)$ for all $n$; then $a^n_{-1}=-a(\var)$, since 
  $\fe^{-1}(a(\var))=\pm a(\var)$. One can achieve this replacing $\ha^n$ 
  by $\hf^{-u_n}(\ha^n)$, $u_n=\max\{ j \ | \ a^n_{-j}=a(\var)\}$. 
  The values $\beta(\ha,\ha^n)$ remain unchanged, by the invariance of 
  basic cocycle. Then $\ha^n$ satisfy the conditions of Lemma \ref{ym}. 
  Therefore, $\beta(\ha,\ha^n)$ converge to a sum (\ref{countsum}) 
  (Corollary \ref{corconv}). This together with the semigroup property 
  (Corollary \ref{semigr}) proves (\ref{countsum}). The set $\be$ is a countable 
  subset in $\rr_{\pm}$, by (\ref{countsum}) and the countability of the set 
  $\Pi_a$. This together with (\ref{closhor}) proves Theorem \ref{tnodense}. 
  \end{proof}

\def\numa{\nu_m(\ha^n)}
\def\numla{\nu_{m+1}(\ha^n)}
\def\nud{\nu_2(\ha^n)}
 
  \subsection{The disk containing the Julia set. Proof of Lemma \ref{jul}}
 We prove Lemma \ref{jul} in the case, when $0<\var<\frac14$. Its proof for $\var<0$ is analogous and 
 will be discussed at the end of the subsection. 
 
 Let $0<\var<\frac14$. 
 We have to show that $J\cap\overline{D_{a(\var)}}=\pm a(\var)$. To do this, we prove that 
 \begin{equation}|\fe(x)|<a(\var) \ \text{for every} \ x\in\partial D_{a(\var)}\setminus\pm a(\var).
 \label{feae}\end{equation}
 This implies that $\fe$ maps the disk $D_{a(\var)}$ to itself, by the maximum 
 principle. Hence, this disk does not intersect the Julia set, by 
 Montel's theorem (\cite{lyu}, p.52). One has 
 $\fe(\overline{D_{a(\var)}}\setminus\pm a(\var))\subset D_{a(\var)}$, by 
 (\ref{feae}). 
 Thus, the set $\overline{D_{a(\var)}}\setminus\pm a(\var)$ is disjoint from 
 $J$, as is 
 $D_{a(\var)}$, by the above inclusion and (\ref{julin}). The points 
 $\pm a(\var)$ belong to the Julia set, since they are mapped to the repelling 
 fixed point $a(\var)$. This proves the second statement of the lemma modulo 
 (\ref{feae}). 
 
 \begin{proof} {\bf of (\ref{feae}).} By definition, $\fe(a(\var))=a^2(\var)+\var=a(\var)$. Hence, for every 
 $\phi\in\rr$ one has
 $$\fe(a(\var)e^{i\phi})=a^2(\var)e^{2i\phi}+\var=e^{2i\phi}(a^2(\var)+\var-\var(1-e^{-2i\phi}))=
 e^{2i\phi}(a(\var)-\var(1-e^{-2i\phi})).$$
 Statement (\ref{feae}) in a equivalent reformulation says that the latter right-hand side has 
 module less than $a(\var)$, whenever $e^{2i\phi}\neq1$. Or equivalently, 
 \begin{equation}
 \Gamma_{\var}=\{\var(1-e^{2i\phi}) \ | \ \phi\in\rr\}\subset 
 D_{a(\var)}(a(\var))\cup0.\label{incga}\end{equation}
 Indeed, the circle $\Gamma_{\var}$ is centered at $\var>0$, 
 has radius $\var$ and is tangent to the boundary 
 $\partial D_{a(\var)}(a(\var))$ at 0. One has  $a(\var)>\var$: 
 the forward 
 orbit of the critical value $\var$ converges to a finite attracting fixed 
 point (\cite{lyu}, p.62), 
 while the orbit of every $x>a(\var)$ tends to infinity. The two latter statements 
 imply (\ref{incga}). 
  This proves (\ref{feae}) and the second statement of  Lemma \ref{jul}. 
 \end{proof}
 
 Let us now consider the case, when $\var<0$. We have to show that 
 $J\subset D_{a(\var)}\cup\pm a(\var)$, or equivalently,  
 $\cc\setminus (D_{a(\var)}\cup\pm a(\var))\subset\cc\setminus J$. To do this, we prove that 
 \begin{equation}|\fe(x)|>a(\var), \ \text{whenever} \ x\in\partial D_{a(\var)}\setminus\pm a(\var). 
 \label{fenew}\end{equation} 
 This is proved analogously to (\ref{feae}) (now $\Gamma_{\var}$ lies outside 
 $D_{a(\var)}(a(\var))$, since $\var<0$). Then 
 \begin{equation}|\fe||_{\cc\setminus(D_{a(\var)}\cup\pm a(\var))}>a(\var).
 \label{fean}\end{equation}
  Indeed, the polynomial $\fe$ has no zeros in the complement
 $\cc\setminus\overline{D_{a(\var)}}$: both its roots $\pm\sqrt{|\var|}$ have module less than $a(\var)$ 
 by (\ref{aedef}). This together with (\ref{fenew}) and the maximum principle 
 applied to the function $\frac1{\fe}(\frac1w)$ implies (\ref{fean}). The 
 complement $\cc\setminus(D_{a(\var)}\cup\pm a(\var))$ is disjoint from 
 the Julia set, by (\ref{fean}) and Montel's theorem, as in the previous case. 
 This proves Lemma \ref{jul}. 
   
    \subsection{Limits of basic cocycles. Proof of Lemma \ref{ym}}

Set $\nu_1=0$. The indices $\nu_m(\ha^n)$ with $m\geq2$ 
are defined below. They are exactly those indices, 
for which $a^n_{-\nu_m}$ is  close to $a(\var)$ (as $n$ is big), and 
$a^n_{-\nu_m-1}$ 
is obtained from it by the inverse branch $\fe^{-1}$ sending $a(\var)$ to $-a(\var)$. 
We show that the number of those $\nu_m$ is finite and uniformly bounded 
and prove Lemma \ref{ym} 
for them. To do this, we choose appropriate $\sigma>0$ small enough, and for 
each $\hy\in\Pi_a$ with $y_{-1}=-a(\var)$ we consider those indices $j$ for which 
$y_{-j}\in D_{\sigma}(a(\var))$ and $y_{-j-1}\notin D_{\sigma}(a(\var))$. 
We show, see (\ref{b>d}), that their number is 
bounded from above by a constant depending only on $\var$, $\sigma$ 
and $\beta(\ha,\hy)$. 

Fix a $\var\in(-\infty,\frac14)\setminus\{-2,0\}$. Fix a $\sigma>0$ and set  
\begin{equation}D'_{\sigma}= \ \text{the connected component containing} \ -a(\var) \ 
\text{of} \ \fe^{-1}(D_{\sigma}(a(\var))).
\label{defks}\end{equation}
We choose $\sigma$ small enough so that 
\begin{equation}
\fe|_{\dse} \ \text{is univalent and} \ \fe(\dse)\supset\overline{\dse},
\label{unidse}
\end{equation}
\begin{equation}\text{the open sets} \ D_{\sigma}(a(\var)), \ D'_{\sigma}, \ 
\fe^{-1}(D'_{\sigma}), \ \fe^{-2}(D'_{\sigma}) 
\ \text{are disjoint. Set}\label{disj}\end{equation} 
\begin{equation}
\delta=\frac12\min\{ |\ln|\fe'(x)|-\ln|\fe'(a(\var))|| \ | \ x\in J(\fe)\setminus(D'_{\sigma}
\cup\dse)\}. \ \text{One has} \ \delta>0.\label{ddelta}\end{equation}
This follows from Corollary \ref{cjul}. By (\ref{unidse}) and (\ref{disj}),
\begin{equation}\fe^{-1}(\dse)=D_{\sigma}'\bigsqcup D_{\sigma}'', \ 
D_{\sigma}''\subset\dse.
 \label{subdse}\end{equation}

Let $\hy\in\Pi_a$, $y_{-1}=-a(\var)$.  Recall that the backward orbit $y_0,y_{-1},\dots$ 
tends to $a(\var)$. Therefore, it leaves $\dse$ only at a finite number of indices, let 
$s=s(\hy)$ denote their number. After each leaving, it returns back either forever, or 
until the next leaving. Let $J(\hy)=\{j_1,\dots,j_s\}$ denote the leaving indices, 
$K(\hy)=\{k_1,\dots,k_s\}$ be the indices of returns:
\begin{equation}J(\hy)=\{ j\in\mathbb N\cup0 \ | \ y_{-j}\in\dse, \ y_{-j-1}\notin\dse\};\label{defJ}
\end{equation}
$$K(\hy)=\{ j\in\mathbb N \ | \ y_{-j+1}\notin\dse, \ y_{-j}\in\dse\}. \text{
 One has}$$
\begin{equation}0=k_0=j_1<k_1\leq j_2<k_2\leq\dots<k_s<j_{s+1}=k_{s+1}=+\infty, 
\ j_r=j_r(\hy), \ k_r=k_r(\hy),\label{defjs}\end{equation}
\begin{equation}y_{-j}\notin\dse, \text{ if and only if } j_r+1\leq j<k_r 
\text{ for some } r\leq s,\label{notind}\end{equation}
\begin{equation} y_{-j}\in D_{\sigma}', \text{ if and only if } j=j_r+1 
\text{ for some } r\leq s,\label{ins'}\end{equation}
\begin{equation} y_{-j}\notin \dse\cup D_{\sigma}', \text{ if and only if } 
j_r+2\leq j<k_r 
\text{ for some } r\leq s,\label{jyj}\end{equation}
\begin{equation} k_r-j_r\geq4 \text{ for each } r\leq s.\label{krjr}\end{equation}
Statements (\ref{defjs}) and (\ref{notind}) follow from definition. Statement 
(\ref{ins'}) follows from definition and (\ref{subdse}). Statement (\ref{jyj}) 
follows from (\ref{notind}) and (\ref{ins'}). One has $y_{-j_r-1}\in D_{\sigma}'$, 
$y_{-j_r-1}, y_{-j_r-2}, y_{-j_r-3}\notin\dse$ for every $r\leq s$, 
by (\ref{disj}) and (\ref{ins'}). Hence, $k_r\geq j_r+4$, by (\ref{notind}). 
This proves (\ref{krjr}). Set 
$$\mathcal D(\hy)=\{ j\in\nn \ | \ y_{-j}\notin\dse\}. 
\text{ By (\ref{notind}),}$$
\begin{equation}\mathcal D(\hy)=\{ j\in\mathbb N \ | \ 
j_r(\hy)+1\leq j<k_{r}(\hy) \text{ for some } 
r\leq s(\hy)\}. \text{ Let}\label{dha}\end{equation}
\begin{equation} d(\hy)=\#(\mathcal D(\hy)) \ \text{ denote its cardinality. One has } 
s(\hy)\leq d(\hy),\label{dy=}\end{equation}
by (\ref{krjr}) and (\ref{dha}).  

The next a priori bound is the main argument in the proof of Lemma \ref{ym}:
  \begin{equation}|\b(\ha,\hy)|>\delta d(\hy)\geq\delta \ \text{for every} \ 
  \hy\in\Pi_a 
  \text{ with } y_{-1}=-a(\var).
  \label{b>d}\end{equation}

\begin{proof} {\bf of (\ref{b>d}).} Each term in (\ref{bxx}) corresponding to an index 
$j$ with $y_{-j}\notin D'_{\sigma}\cup \dse$ 
has module no less than $2\delta$, by the definition of $\delta$. These 
are exactly the indices from (\ref{jyj}), and their number equals 
$\sum_{r=1}^s(k_r-j_r-2)$. 
All the  terms in (\ref{bxx}) are nonzero  and have the same sign, as 
$\b(\ha,\hy)$, except for zero terms with $j=0,1$: $y_{0,-1}=\pm a(\var)$. This 
follows from Corollary \ref{corsign}. 
This together with the previous statements and (\ref{krjr}) implies that  
$$\b(\ha,\hy)\geq2\delta\sum_{r=1}^s(k_r-j_r-2)>\delta\sum_{r=1}^s(k_r-j_r-1)=
\delta d(\hy).$$
\end{proof}

The numbers $d(\ha^n)$ and  $s(\ha^n)$ are uniformly bounded from above, 
by (\ref{dy=}), (\ref{b>d}) and the convergence of the values $\beta(\ha,\ha^n)$. 
Now passing to a subsequence, one can achieve  that 
\begin{equation} s=s(\ha^n), \ d=d(\ha^n) \  \text{are independent on} \ n,\label{sn}\end{equation}
\begin{equation}\text{for every fixed} \ j\in\nn\cup0 \ \text{and} \ r\leq s
\ \text{the sequence} \ a^n_{-j_r(\ha^n)-j} \ \text{converges, as} \ n\to\infty.\label{hanconv}
\end{equation}
Take those values $r$ (denote $l$ their number) for which 
\begin{equation}a^n_{-j_r(\ha^n)}\to a(\var), \ \text{as} \ n\to\infty; \ 
\text{denote these} \ r \ \text{by} \ 
1=r(1)<\dots<r(l). \text{ Set } \nu_{l+1}=+\infty,\label{rm}\end{equation}
 \begin{equation}\nu_m(\ha^n)=j_{r(m)}(\ha^n), \ y_{-j}(m)=
 \lim_{n\to\infty}a^n_{-\nu_m(\hat a^n)-j} \ 
 \text{for every} \ m=1,\dots,l,  \ j\in\nn\cup0,
 \label{defnu}\end{equation} 
 $$\hy(m)=(y_0(m),y_{-1}(m),\dots).$$
  Statements 
 (\ref{convnu})-(\ref{defam}) follow from (\ref{rm}) and (\ref{defnu}). 

\begin{proof} {\bf of (\ref{pia}).} 
One has $y_0(m)=a(\var)$, by (\ref{rm}),   
$y_{-j}(m)\to a(\var)$, as $j\to+\infty$ (hence, $\hy(m)\in L(\ha)$). 
Indeed, the converse would imply that 
there  exists a $\sigma'>0$ such that there is an infinite number of indices $j$ 
for which $y_{-j}(m)\notin\overline{D_{\sigma'}(a(\var))}$. This together with the 
convergence $a^n_{-\nu_m(\ha^n)-j}\to y_{-j}(m)$ would 
imply that there are $\ha^n$ with arbitrarily large number of indices $j$ 
for which $ a^n_{-j}\notin D_{\sigma'}(a(\var))$. Hence, the number $d(\ha^n)$ 
constructed above with $\sigma$ replaced by $\min\{\sigma,\sigma'\}$  is not uniformly 
bounded, - a contradiction to (\ref{b>d}) and the convergence of 
$\beta(\ha,\ha^n)$. Thus, $\hy(m)\in L(\ha)$, 
$y_0(m)=a(\var)$. One has $\hy(m)\neq\ha$, since $y_{-1}(m)=-a(\var)$. 
Indeed, $y_{-1}(m)\in\{\pm a(\var)\}=\fe^{-1}(a(\var))$, 
since $y_0(m)=a(\var)$. By definition, the point $y_{-1}(m)$ is the limit of 
$a^n_{-1-\nu_m(\ha^n)}\notin\dse$, see (\ref{defJ}) and 
(\ref{defnu}). Hence, $y_{-1}(m)=-a(\var)$. 
The projection $\pi_0|_{L(\ha)}$ has nonzero 
derivative at $\hy(m)$ (hence, $\hy(m)\in\Pi_a$). 
Indeed, otherwise some $y_{-j}(m)$ equals 0: the critical point of 
$\fe$. Then there are subsequences of indices $n_u$ and $q_u$ such that 
$a^{n_u}_{-q_u}\to0$, as $u\to\infty$. Therefore,  the  terms in (\ref{bxx}) 
corresponding to the points $y_{-j}=a^{n_u}_{-q_u}$ 
tend to infinity. Hence, $\b(\ha,\ha^{n_u})\to\infty$ 
(Corollary \ref{corsign}), - a  contradiction. Statement (\ref{pia}) is proved. 
\end{proof}

\begin{proof} {\bf of (\ref{sigmam}).} Fix a $m=1,\dots,l$. 
Everywhere in the proof, whenever the contrary is not specified, 
$\fe^{-1}$ denotes the inverse branch $\dse\to\dse$  fixing $a(\var)$: its iterates 
converge to $a(\var)$ uniformly on $\dse$, by (\ref{unidse}).  
Recall that $a^n_{-\nu_k(\ha^n)}\in\dse$ for all $k=1,\dots,l$ and 
$n\in\mathbb N$, by (\ref{defJ}). 
Consider the following auxiliary backward orbits. If $m<l$, we set    
\begin{equation}
\hb^n=( a^n_{-\nu_m(\ha^n)}, a^n_{-\nu_m(\ha^n)-1},\dots, 
a^n_{-\nu_{m+1}(\ha^n)},
\fe^{-1}( a^n_{-\nu_{m+1}(\ha^n)}),\fe^{-2}( a^n_{-\nu_{m+1}(\ha^n)}),\dots),
\label{bnm}
\end{equation}
$$\hc^n=\hat\fe^{\nu_m(\ha^n)-\nu_{m+1}(\ha^n)}(\hb^n)=
( a^n_{-\nu_{m+1}(\ha^n)},
\fe^{-1}( a^n_{-\nu_{m+1}(\ha^n)}),\dots). \ \text{If} \ m=l, \ \text{we set}$$
$$\hb^n=\hat\fe^{-\nu_l(\ha^n)}(\ha^n)=( a^n_{-\nu_l(\ha^n)}, 
a^n_{-\nu_l(\ha^n)-1},\dots), \ 
\hc^n=\ha.$$
By construction, $\hb^n,\hc^n\in L(\ha)$.  We use the formula 
\begin{equation}\Sigma_m=\b_{\hb^n}-\b_{\hc^n}.\label{sig=-}\end{equation}
This follows from definition, see (\ref{sigmam}), and formula (\ref{bxx'}) applied to the following pairs of orbits: 
$\hx=\ha$, $\hy=\hb^n$; \ $\hx=\ha$, $\hy=\hc^n$. We show that  
\begin{equation}
\b_{\hb^n}\to\b_{\hy(m)}, \ \b_{\hc^n}\to\b_{\ha}, \ \text{as} \ n\to\infty.
\label{bhbn}
\end{equation}
This together with (\ref{sig=-}) and the equality $\b(\ha,\hy(m))=\b_{\hy(m)}-\b_{\ha}$ 
implies (\ref{sigmam}).

For every $n\in\mathbb N$ one has $\hc^n\in L(\ha,\dse)$, and 
$c^n_0=a^n_{-\nu_{m+1}(\ha^n)}\to a(\var)$, as $n\to\infty$ 
($\hc^n=\ha$, if $m=l$), by construction and (\ref{rm}). Hence, $\hc^n\to\ha$,    
$\b_{\hc^n}\to\b_{\ha}$, as $n\to\infty$. To prove that $\b_{\hb^n}\to\b_{\hy(m)}$, we fix a neighborhood 
$V=V(a(\var))\subset\oc$ such that the local 
leaf $L(\hy(m),V)$ is univalent over $V$. We show that 
\begin{equation}\hb^n\in L(\hy(m),V) \text{ for every } n \text{ large enough.}
\label{blocl}\end{equation} 
\begin{proof} {\bf of (\ref{blocl}).} If $m<l$, we set
$$\mu(\ha^n)=k_{r(m+1)-1}(\ha^n), \ 
\text{see (\ref{defjs}). If} \ m=l, \ \text{we set} \ 
\mu(\ha^n)=k_s(\ha^n).$$ 
One has $\nu_m(\ha^n)=j_{r(m)}(\ha^n)<\mu(\ha^n)\leq\nu_{m+1}(\ha^n)=
j_{r(m+1)}(\ha^n)$, by (\ref{defjs}), 
\begin{equation}a^n_{-\mu(\ha^n)},\dots,a^n_{-\nu_{m+1}(\ha^n)}\in\dse, 
\label{indse}\end{equation}
by (\ref{notind}). Hence, by (\ref{bnm}), 
\begin{equation}
b^n_{-j}\in\dse, \  b^n_{-j-1}=\fe^{-1}( b^n_{-j}) \ \text{for every} \ j\geq 
q(n)=q_m(n)=\mu(\ha^n)-\nu_m(\ha^n).\label{qnm}\end{equation}

{\bf Claim 1.} {\it The above numbers $q(n)$ are uniformly bounded from above.} 

\begin{proof} There exists a $\sigma'>0$ such that  
\begin{equation}b^n_{-j}=a^n_{-j-\nu_m(\ha^n)}\notin D_{\sigma'}(a(\var)) 
\text{ for every } n\in\mathbb N \text{ and } 0<j<q(n). \label{sigma'b}
\end{equation}
Indeed, the  
converse would mean that there exist sequences of indices $n_u,p_u\in\mathbb N$, 
$\nu_m(\ha^{n_u})<p_u<\mu(\ha^{n_u})$, such that $a^{n_u}_{-p_u}\to a(\var)$, 
as $u\to\infty$. Note that the latter lower and upper bounds for $p_u$ are the 
moments of leaving $\dse$ and returning to $\dse$ of the backward orbit $\ha^{n_u}$. 
Without loss of generality we consider that the holomorphic branch $\fe^{-1}:\dse
\to\dse$ contracts the distances. Then for every $\hy\in\Pi_a$ with 
$y_{-1}=-a(\var)$ 
the local minima for $dist(y_{-j},a(\var))$ are achieved at $j\in J(\hy)$. Hence, 
after passing to a subsequence, the above $p_u$ may be chosen equal to 
$j_r(\ha^{n_u})$, $r(m)<r<r(m+1)$ ($r>r(m)$, 
if $m=l$): $a^{n_u}_{-j_r(\ha^{n_u})}\to a(\var)$, as $n\to\infty$, - a contradiction 
to the definition of the numbers $r(m)$, see (\ref{rm}). 
This proves (\ref{sigma'b}). It implies that the number $q(n)-1$ is no greater 
than the number $d(\ha^n)$ from (\ref{dy=}) defined with $\sigma$ 
replaced by $\sigma'$. This together with (\ref{b>d}) proves the claim.
\end{proof}
 
Without loss of generality we assume that $q=q(n)$ is independent on $n$, 
passing to a subsequence. 
Let $V$ be the neighborhood of $a(\var)$ from (\ref{blocl}). 
For the proof of (\ref{blocl}) we 
have to show that for every $n$ large enough and each $j\in\nn$ the inverse 
branch $\fe^{-j}:y_0(m)\mapsto y_{-j}(m)$ (which is single-valued on $V$ by 
construction) 
sends $ b_0^n$ to $ b_{-j}^n$. Indeed, for every $j\in\mathbb N\cup0$ one has 
$b^n_{-j}\to y_{-j}(m)$, as $n\to\infty$, by (\ref{defnu}) and (\ref{bnm}). 
This implies the latter inverse branch statement for each $n$ large enough and 
every $j\leq q$. For each $j>q$ the attractive inverse branch 
$\fe^{-1}:\dse\to\dse$ sends $b^n_{-j}$ to $b^n_{-j-1}$ for every 
$n\in\mathbb N$, by (\ref{bnm}), (\ref{indse}) and (\ref{subdse}). 
The same holds true for $b^n_{-j}$ replaced by $y_{-j}(m)$, passing to the limit. 
This proves (\ref{blocl}). 
\end{proof}

One has $\hb^n\to\hy(m)$, $\b_{\hb^n}\to\b_{\hy(m)}$, as $n\to\infty$, 
by (\ref{blocl}) and since $b^n_0=a^n_{-\nu_m(\ha^n)}\to y_0(m)=a(\var)$. 
This  proves (\ref{bhbn}) and (\ref{sigmam}). The proof of Lemma \ref{ym} is complete.
\end{proof}

\subsection{Accumulation set of a horosphere. Proof of (\ref{closhor})} 
For the proof of (\ref{closhor}) we have to prove the  two following statements:

\begin{equation}\text{the accumulation set of the horosphere} \ S=S_{\ha,0} \ 
\text{is contained in} \ SB_\var;\label{accinc}\end{equation}
\begin{equation} \text{the horosphere} \ S \ \text{accumulates at least to the whole set} \ SB_\var.
\label{acc=}\end{equation}

Recall that we measure heights of horospheres with respect to the Euclidean 
metric of $\cc$ lifted to the leaves of $\cafe$. 

\begin{proof} {\bf of (\ref{acc=}).} For every $\hy\in\Pi_a$ and $n\in\mathbb N$ set
$\a^n(\hy)=(\hf^n(\hy),\beta(\ha,\hy))\in\chfe$. One has $\a^n(\hy)\in S$, 
since $(\hf^n(\hy),\beta(\ha,\hf^n(\hy))\in S$, see (\ref{monodr}), and 
$\beta(\ha,\hy)=\beta(\ha,\hf^n(\hy))$ 
(the invariance of basic cocycle); $\a^n(\hy)\to\a(\hy)=(\ha,\beta(\ha,\hy))$, as 
$n\to+\infty$  
($\hf^n(\hy)\to\ha$, by Corollary \ref{cconva}). Hence, $S$ accumulates at least to 
the horospheres through all the points $\a(\hy)$, $\hy\in\Pi_a$. 
The set $SB_\var$ is the closure of the union of the latter horospheres. 
This proves (\ref{acc=}).
\end{proof}

As it is shown below, statement (\ref{accinc}) is implied by the following 
proposition. 
\begin{proposition} \label{propohm} 
Let $\ha^n$, $l$, $\nu_m$, $\hy(m)$ be the same, as in Lemma \ref{ym}. Then 
 $\ha^n\to\hy(1)$ in $\cafe$, as $n\to\infty$. 
\end{proposition}

\begin{corollary} \label{cpropohm} For every 
sequence of distinct points $(\ha^n,h_n)\in S$, $\ha^n\in\Pi_a$, converging in 
$\chfe$ to  a point 
$(\hy,h)\in\chfe$ with $(\pi_0|_{L(\hy)})'(\hy)\neq0$, $h\in\rr$, 
the latter limit  is contained in $SB_\var$. 
\end{corollary}

In the proof of Proposition \ref{propohm} and Corollary \ref{cpropohm} we use 
the following 

{\bf Claim 2.} {\it In the conditions of Lemma \ref{ym} for every $n$ large enough and 
each $j\leq\nu_2(\ha^n)$ one 
has $a^n_{-j}=y_{-j}(1)$. For every $m=1,\dots,l$ one has $a^n_{-\nu_m(\ha^n)-j}
\to y_{-j}(m)$, as $n\to\infty$, uniformly in $0\leq j\leq\numla-\numa$.}

\begin{proof} Let $q(n)=q_m(n)$ be the numbers from (\ref{qnm}). They are uniformly 
bounded in $m$ and $n$, by Claim 1 at the same place. Let $q=\max_{m,n}q_m(n)$. For every 
$q\leq j\leq\numla-\numa$ one has $a^n_{-\numa-j}\in\dse$, by (\ref{indse}). 
Hence, for each $q<j\leq\numla-\numa$ the point $a^n_{-\numa-j}$ is obtained 
from $a^n_{-\numa-j+1}$ by applying the attractive inverse 
branch $\fe^{-1}:\dse\to\dse$, by (\ref{subdse}). The same branch transforms the limit 
$y_{-j+1}(m)$ of the latter to the limit $y_{-j}(m)$ of the former, see 
(\ref{hanconv}). Therefore, the convergence $a^n_{-\nu_m(\ha^n)-j}
\to y_{-j}(m)$, which is uniform in $j\leq q$, is also uniform in all  
$j\leq\numla-\numa$. 
This proves the uniform convergence statement of the claim. The first statement 
of the claim holds true by the uniform convergence and since 
$a^n_0=y_0(m)=a(\var)$. 
\end{proof}

\begin{proof} {\bf of Corollary \ref{cpropohm}.} Let $\a^n=(\ha^n,h_n)$ be the same, 
as in the corollary. Then $h_n=\beta(\ha,\ha^n)$, since $\a^n\in S$ by assumption 
and by (\ref{monodr}). Set 
$$u_n=\max\{ j \ | \ a^n_{-j}=a(\var)\}, \ \hb^n=\hf^{-u_n}(\ha^n):$$
\begin{equation}b^n_{-1}=-a(\var), \ \beta(\ha,\hb^n)=\beta(\ha,\ha^n)=h_n\to h, 
\text{ as } n\to\infty,\label{bbe}\end{equation} 
by the invariance of basic cocycle. 
Thus, the sequence $\hb^n$ satisfies the conditions 
of Lemma \ref{ym}. Let $l\in\mathbb N$ and $\hy(m)\in\Pi_a$ be the same as in 
the lemma. Then $\hb^n\to\hy(1)$ in $\cafe$ (Proposition \ref{propohm}), 
\begin{equation}(\hb^n,h_n)\to 
(\hy(1),h) \text{ in }  \chfe, \   h=\sum_{m=1}^l\beta(\ha,\hy(m))\label{hbhb}\end{equation}
 (Corollary \ref{corconv}). One has $(\hb^n, h_n)=(\hb^n,\beta(\ha,\hb^n))\in S$, by 
(\ref{bbe})  and (\ref{monodr}).  Let us consider the three following possible cases:

Case 1): $u_n=0$ for all $n$. Then $\ha^n=\hb^n$. One has $l\geq2$. Indeed, 
otherwise, if $l=1$, then $\hb^n=\hy(1)$ for every $n$ large enough 
(Claim 2: here $\nu_2=\nu_{l+1}=+\infty$). But $(\hb^n,h_n)\in S$ are distinct points, by assumption. Hence, 
$\hb^n$ are also distinct, since the restriction to 
each horosphere of the $\pi_{hyp}$- projection is injective. This  
contradicts the equality $\hb^n=\hy(1)$. Thus, $l\geq 2$, $\hy=\hy(1)$, 
\begin{equation}S_{\hy,h}=S_{\ha,h'} , \ 
h'=h-\beta(\ha,\hy)=\sum_{m=2}^l\beta(\ha,\hy(m)),\label{dddd}\end{equation} 
 by (\ref{hbhb}) and (\ref{monodr}). The latter horosphere is contained in $SB_\var$ by 
(\ref{countsum}). This proves the statement of the corollary.

Case 2): $u_n$ is bounded. Without loss of generality we consider that 
$u_n=const=u$, passing to a subsequence. Then $l\geq2$, as in the above case, 
and $\hy=\hf^u(\hy(1))$. Statement (\ref{dddd}) holds true again, since 
$$h-\beta(\ha,\hy)=\sum_{m=1}^l\beta(\ha,\hf^u(\hy(m)))-\beta(\ha,\hf^u(\hy(1))
=h',$$
by (\ref{hbhb}) and the equalities $\beta(\ha,\hf^u(\hy(m)))=\beta(\ha,\hy(m))$ 
(the invariance of basic cocycle). This proves the corollary.

Case 3): $u_n\to\infty$ (after passing to a subsequence). Then 
\begin{equation} \ha^n=\hf^{u_n}(\hb^n)\to\ha \text{ in } \cafe, 
\text{ as } n\to\infty.\label{convbn}\end{equation}
Indeed, let $V\subset\cc$ 
be a neighborhood of $a(\var)$ over which the local leaf $L(\hy(1),V)$ is 
univalent, $U$ be a smaller neighborhood: $\overline U\subset V$. 
Then the local leaves $L(\hb^n,U)$ are univalent over $U$, whenever $n$ is 
large enough, by the convergence $\hb^n\to\hy(1)$ in $\cafe$ and Proposition 
\ref{conva}. This together with Corollary \ref{cconva} implies (\ref{convbn}). 
Thus, $(\ha^n,h_n)\to(\ha,h)$, $h_n\in B_{\var}$, and hence, 
$h\in B_{\var}$, since $B_{\var}$ is closed. Therefore, 
$(\ha,h)\in SB_{\var}$. The proof of the corollary is complete. 
\end{proof}

\begin{proof} {\bf of (\ref{accinc}) modulo Proposition \ref{propohm}.} 
Let the horosphere $S$ accumulate to some horosphere $S'\subset\chfe$, 
$L'\subset\cafe$ be the corresponding affine leaf. The leaf $L'$ is contained in 
$\cafe'$, and it is dense there (Proposition \ref{minim}). Hence, there exists a point 
$\hy\in L'$ 
such that $\pi_0(\hy)=a(\var)$ and $(\pi_0|_{L'})'(\hy)\neq0$, since 
$\Pi_a\neq\emptyset$. Let $h$ be the height of the 
horosphere 
$S'$ over $\hy$. Then $(\hy,h)$ is a limit of a sequence of distinct points 
$(\ha^n,h_n)\in S$, 
$\ha^n\in\Pi_a$. Hence, $(\hy,h)\in SB_\var$ and $S'\subset SB_\var$, by 
Corollary \ref{cpropohm}. This proves (\ref{accinc}).
\end{proof}
 
\begin{proof} {\bf of Proposition \ref{propohm}.} One has $\ha^n\to\hy(1)$ in 
$\mathcal N_{\fe}$, by (\ref{hanconv}). Let us prove the convergence in $\cafe$. 
Fix a $N\in\mathbb N$ and a neighborhood 
$V$ of $y_{-N}(1)$ such that the local leaf $L(\hat\fe^{-N}(\hy(1)),V)$ is 
univalent over $V$, and its arbitrary smaller neighborhood $V'$, 
$\overline{V'}\subset V$. Let us show that the local leaf 
$L(\hat\fe^{-N}(\ha^n),V')$ is also  univalent over $V'$, whenever 
$n$ is large enough. This together with Proposition \ref{conva} proves 
Proposition \ref{propohm}. 

Fix a neighborhood $W=W(a(\var))\subset\cc$  such that for every 
$m=2,\dots,l$ the local leaf $L(\hy(m),W)$ is  univalent over $W$.  
Fix an arbitrary smaller neighborhood $U=U(a(\var))$, $\overline{U}\subset W$. 
For every $k\in\nn\cup0$ and each $\hx\in\cafe^l$ set
$$f_{\var,\hx}^{-k}=\text{ the germ of the inverse branch } \fe^{-k} \text{ sending } 
x_0 \text{ to } x_{-k}. \text{ Set}$$
$$g^{-k}_1=f^{-k}_{\var,\hf^{-N}(\hy(1))}, \ g_m^{-k}=f_{\var,\hy(m)}^{-k} 
\text{ for every } m\geq2,$$
$$g_{N,n}^{-k}=f_{\var,\hf^{-N}(\ha^n)}^{-k}, \ 
\phi_{m,n}^{-k}= f_{\var,\hf^{-\numa}(\ha^n)}^{-k}\text{ for every } 
n\in\mathbb N, \ m=1,\dots,l. 
\text{ One has:}$$
\begin{equation}
g_1^{-k} \text{ is holomorphic on } V \text{ and } g_1^{-k}\to a(\var) 
\text{ uniformly on } \overline{V'}, \text{ as } k\to\infty,\label{hol1}\end{equation}
\begin{equation}
g_m^{-k}  \text{ is holomorphic on } W \text{ and } g_m^{-k}\to a(\var) 
\text{ uniformly on } \overline{U}, \text{ for every } m\geq2,\label{holm}\end{equation}
 by the univalence of the local leaves $L(\hat\fe^{-N}(\hy(1)),V)$,  
 $L(\hy(m),W)$ and the Shrinking Lemma. 
 
 \medskip
 
 {\bf Claim 3.} {\it For every $n$ large enough and $k\in\mathbb N$ 
 the function $g_{N,n}^{-k}$ is holomorphic on $V'$.}
 
 \begin{proof} For every $n$ large enough  one has 
 \begin{equation}a^n_{-N-k}=y_{-N-k}(1), \text{ thus, } g_{N,n}^{-k}=g_1^{-k},  
 \text{ for every } k\leq\nud-N,\label{AA}
 \end{equation} 
 \begin{equation}g_1^{N-\nud}(\overline{V'})\subset U, \  
 g_m^{\nu_m(\ha^n)-\nu_{m+1}(\ha^n)}(\overline U)\subset U \text{ for every } 
 m\geq2,\label{feincl}\end{equation}
 \begin{equation} g_m^{-k}=\phi_{m,n}^{-k} \text{ for every } m\geq2 
 \text{ and } k\leq\numla-\numa.\label{numl}\end{equation}
 Statement (\ref{AA}) follows from the first statement of Claim 2. Statement 
 (\ref{feincl}) follows from 
 (\ref{hol1}), (\ref{holm}) and (\ref{diffinf}). Statement (\ref{numl}) follows from 
 the second statement of Claim 2. Fix an arbitrary $n$ satisfying 
 (\ref{AA})-(\ref{numl}). Let us prove the 
 holomorphicity of $g_{N,n}^{-k}$ on $V'$. We prove this for every $k>\nu_l(\ha^n)$: 
 this implies the holomorphicity of $g_{N,n}^{-k'}=f^{k-k'}\circ g_{N,n}^{-k}$ for 
 smaller $k'$. One has $g_{N,n}^{N-\nud}=g_1^{N-\nud}$, by (\ref{AA}). Hence, 
 $$g_{N,n}^{-k}=\phi_{l,n}^{\nu_l(\ha^n)-k}\circ 
 \phi_{l-1,n}^{\nu_{l-1}(\ha^n)-\nu_l(\ha^n)}\circ
 \dots\circ\phi_{2,n}^{\nu_2(\ha^n)-\nu_3(\ha^n)}\circ 
 g_{N,n}^{N-\nud}$$
 $$=g_l^{\nu_l(\ha^n)-k}\circ g_{l-1}^{\nu_{l-1}(\ha^n)-\nu_l(\ha^n)}\circ
 \dots\circ g_2^{\nu_2(\ha^n)-\nu_3(\ha^n)}\circ g_1^{N-\nud},$$
 and the latter composition is well-defined and holomorphic on $V'$, by (\ref{feincl}) 
 and (\ref{numl}).  This together with the previous discussion proves Claim 3.
 \end{proof}

Claim 3 implies the univalence of the local leaf $L(\hat\fe^{-N}(\ha^n),V')$ 
over $V'$ for every $n$ large enough. Together with Proposition \ref{conva}, this 
proves Proposition \ref{propohm}.
\end{proof}

 \subsection{Perturbations. Proof of the Addendum to Theorem \ref{tnodense}}
The proof of the Addendum to Theorem \ref{tnodense} is analogous to that of 
the theorem itself with  minor modifications listed below. Most of them 
concern the proof of inequality (\ref{b>d}).  We fix a $\var_0<\frac14$, 
$\var_0\neq0$ (that may be equal to $-2$). Let $\sigma$ and $\delta$ be the 
numbers from (\ref{unidse})-(\ref{ddelta}), both corresponding to 
$\var=\var_0$. We show that there exists a $\Delta>0$ such that 
(\ref{b>d}) holds true  for every complex value  $\var\neq-2$, 
$|\var-\var_0|<\Delta$, with $\delta$ replaced by $\frac{\delta}2$. 
To do this, we use the three following well-known facts:

- the closure of a (super-)attracting basin of a rational function of a given degree 
depends lower-semicontinuously on its coefficients in the Hausdorff topology;

- the Julia set of a complex polynomial is the boundary of the super-attracting 
basin of infinity; 

-  there is a complex neighborhood of the interval 
$0<\var<\frac14$ where the Julia set $J(\fe)$ 
is a Jordan curve separating the super-attracting basin of infinity from the 
attracting basin  of a finite attracting fixed point.

For every 
given $\var_0<0$ ($0<\var_0<\frac14$) and each $\sigma'>0$ there exists a 
$\Delta>0$ such that 
for every $\var\in\cc$, $|\var-\var_0|<\Delta$, the corresponding inclusion from 
Lemma \ref{jul} holds true with $J(\fe)$ replaced by 
$$J_{\var,\sigma'}=J(\fe)\setminus(D_{\sigma'}(a(\var))\cup D_{\sigma'}'):$$ 
\begin{equation}J_{\var,\sigma'}\subset D_{a(\var)} \ (\text{respectively}, \ 
J_{\var,\sigma'}
\subset\cc\setminus\overline{D_{a(\var)}}),\label{newinc}
\end{equation}
where $D'_{\sigma'}$ is the corresponding domain (\ref{defks}). 
This follows from Lemma \ref{jul} and the above semicontinuity statement  
applied to the basin of infinity (respectively, the basin of a finite
 attracting fixed point). 
Let for simplicity $0<\var_0<\frac14$: then for every $\hy\in\Pi_a$ all the terms in the corresponding formula (\ref{bxx}) are positive. 
If  $\var\in D_{\Delta}(\var_0)$ is not real, then, in general, some terms 
in (\ref{bxx}) may be negative. But this may happen only to 
 those corresponding to $y_{-j}\in D'_{\sigma'}\cup D_{\sigma'}(a(\var))$. This follows  
from (\ref{newinc}), as in the proof of Corollary \ref{cjul}. The other terms are no less 
than $\delta$, see (\ref{ddelta}),  
if $\var=\var_0$, and greater than $\frac{2\delta}3$, whenever 
$\var$ is complex and $\var-\var_0$ is small enough. Elementary estimates
show that if $\var-\var_0$ is small enough, then the sum of the latter terms 
strongly dominates that of the former ones so that inequality (\ref{b>d}) holds 
with $\delta$ replaced by $\frac{\delta}2$. The rest of the proof of Theorem 
\ref{tnodense} applies with obvious changes.

 \section{Acknowledgments} I wish to thank M.Yu.Lyubich for attracting my attention 
 to the problem and helpful discussions. 
 I wish to thank  also  D.Saric and C.Cabrera for helpful discussions. The main results 
 of the paper were partially obtained while I was visiting the 
 Institute for Mathematical Sciences, SUNY, Stony Brook (USA) and Fields Institute (Toronto, Canada).  I wish to thank both Institutes for hospitality and support.

\end{document}